\documentclass{article}

\usepackage{microtype}
\usepackage{graphicx}
\usepackage{subfigure}
\usepackage{booktabs} 
\usepackage{makecell}

\usepackage{hyperref}

\usepackage[accepted]{icml2024-1}

\usepackage{amsmath}
\usepackage{amssymb}
\usepackage{mathtools}
\usepackage{amsthm}

\usepackage[capitalize,noabbrev]{cleveref}

\theoremstyle{plain}
\newtheorem{theorem}{Theorem}[section]
\newtheorem{proposition}[theorem]{Proposition}
\newtheorem{lemma}[theorem]{Lemma}

\theoremstyle{definition}
\newtheorem{definition}[theorem]{Definition}
\newtheorem{assumption}[theorem]{Assumption}
\theoremstyle{remark}
\newtheorem{remark}[theorem]{Remark}
\newcommand{\close}{\hfill{\footnotesize$\Diamond$}}

\usepackage[textsize=tiny]{todonotes}

\usepackage{acronym}

\newcommand{\dd}{\:d}

\newcommand{\eps}{\varepsilon}

\newcommand{\dif}{\dd}

\DeclareMathOperator*{\argmin}{argmin}
\DeclareMathOperator*{\Argmin}{Argmin}
\DeclarePairedDelimiter{\ceil}{\lceil}{\rceil}

\DeclareMathOperator{\diag}{diag}

\DeclareMathOperator{\dom}{dom}

\DeclareMathOperator{\Int}{int}

\DeclareMathOperator{\image}{im}

\newcommand{\ct}{\mathtt{t}}

\newcommand{\bA}{{\mathbf A}}
\newcommand{\BB}{{\mathbf B}}
\newcommand{\bH}{\mathbf{H}}
\newcommand{\bI}{\mathbf{I}}
\newcommand{\bJ}{\mathbf{J}}

\newcommand{\gbold}{\mathbf{g}}

\newcommand{\bZ}{\mathbf{Z}}

\renewcommand{\iff}{\Leftrightarrow}

\newcommand{\eqdef}{\triangleq}

\newcommand{\setE}{\mathbb{E}}
\newcommand{\scrF}{\mathcal{F}}

\newcommand{\scrH}{\mathcal{H}}

\newcommand{\setK}{\mathsf{K}}

\newcommand{\setL}{\mathsf{L}}

\newcommand{\scrT}{\mathcal{T}}

\newcommand{\scrW}{\mathcal{W}}

\newcommand{\setX}{\mathsf{X}}

\newcommand{\0}{\mathbf{0}}
\newcommand{\1}{\mathbf{1}}

\newcommand{\R}{\mathbb{R}}

\newcommand{\K}{\mathbb{K}}

\DeclareMathOperator{\NC}{\mathsf{NC}}

\newcommand{\feas}{\mathsf{X}}							

\DeclarePairedDelimiter{\abs}{\lvert}{\rvert}

\DeclarePairedDelimiter{\inner}{\langle}{\rangle}
\DeclarePairedDelimiter{\norm}{\lVert}{\rVert}

\newacro{SC}[SC]{self-concordant}
\newacro{SCB}[SCB]{self-concordant barrier}
\newacro{SSB}[SSB]{self-scaled barrier}
\newacro{FOM}{First-order method}

\DeclareMathOperator{\AHBA}{\mathbf{FOABM}}
\DeclareMathOperator{\SAHBA}{\mathbf{SOABM}}

\def\rev#1{{\color{black}#1}} 

\icmltitlerunning{Barrier Algorithms for Constrained Non-Convex Optimization}

\begin{document}

\twocolumn[
\icmltitle{Barrier Algorithms for Constrained Non-Convex Optimization}




\begin{icmlauthorlist}
\icmlauthor{Pavel Dvurechensky}{wias}
\icmlauthor{Mathias Staudigl}{uman}
\end{icmlauthorlist}

\icmlaffiliation{wias}{Weierstrass Institute
for Applied Analysis and Stochastics, Berlin, Germany}
\icmlaffiliation{uman}{University of Mannheim, Mannheim, Germany}

\icmlcorrespondingauthor{Pavel Dvurechensky}{pavel.dvurechensky@wias-berlin.de}

\icmlkeywords{Machine Learning, ICML}

\vskip 0.3in
]



\printAffiliationsAndNotice{\icmlEqualContribution} 

\begin{abstract}
In this paper we theoretically show that interior-point methods based on self-concordant barriers possess favorable global complexity beyond their standard application area of convex optimization. To do that we propose first- and second-order methods for non-convex optimization problems with general convex set constraints and linear constraints. Our methods attain a suitably defined class of approximate first- or second-order KKT points with the worst-case iteration complexity similar to unconstrained problems, namely  $O(\eps^{-2})$ (first-order) and $O(\eps^{-3/2})$ (second-order), respectively.
\end{abstract}

\section{Introduction}
\label{sec:intro}
Interior-point methods are a universal and very powerful tool for convex optimization \cite{NesIPM94,BoyVan04} that allows obtaining favorable global complexity guarantees for a variety of problems with many applications. Much less is known about the complexity guarantees for such methods in the non-convex world, especially important in machine learning applications such as training neural networks. This paper aims to fill this gap in the theoretical analysis of interior-point methods in application to optimization with non-convex objectives. 

Let $\setE$ be a finite-dimensional vector space with inner product $\inner{\cdot,\cdot}$ and  Euclidean norm $\norm{\cdot}$. Our goal is to solve constrained optimization problems of the form 
\begin{equation}\label{eq:Opt}\tag{Opt}
 \min_{x} f(x)\quad \text{s.t.: } \bA x=b,\; x\in\bar{\setK}.
\end{equation}
Our main  assumption   is as follows:
\begin{assumption}\label{ass:1}
\begin{enumerate}
\item $\bar{\setK}\subset\setE$ is a closed convex  set with nonempty relative interior $\setK$;
\item $\bA:\setE\to\R^{m}$ is a linear operator assigning each element $x\in\setE$ to a vector in $\R^m$ and having full rank, i.e., $\image(\bA)=\R^{m}$, $b \in \R^{m}$; 
\item The feasible set $\bar{\setX}=\bar{\setK}\cap\setL$ with $\setL=\{x\in\setE\vert\bA x=b\}$ has nonempty relative interior denoted by $\setX=\setK\cap\setL$;
\item $f:\setE\to\R$ is possibly \textit{non-convex}, continuous on $\bar{\setX}$ and continuously differentiable on $\setX$;
\item Problem \eqref{eq:Opt} admits a global solution. We let $f_{\min}(\setX)=\min\{f(x)\vert x\in\bar{\setX}\}$.
\end{enumerate}
\end{assumption}

As a main tool for developing our algorithms, we use a \emph{self-concordant barrier} (SCB) $h(x)$ for the set $\bar{\setK}$, see \cite{NesNem94} and Definition \ref{def:LHSCB}. Using the barrier $h$, our algorithms are designed to reduce the \emph{potential function}
\begin{equation}\label{eq:potential}
F_{\mu}(x)=f(x)+\mu h(x),
\end{equation}
where $\mu>0$ is a (typically) small penalty parameter. 

This approach has the following advantages compared to other approaches for solving \eqref{eq:Opt}.
\begin{enumerate}
    \item Unlike proximal-operator- or projection-based approaches \cite{GhaLan16,bogolubsky2016learning,cartis2019optimality,cartis2012complexity,CurRobSam17,CarGouToi12,birgin2017regularization,CarGouToi18,cartis2019optimality} our approach does not require evaluation of a costly projection onto the feasible set $\bar{\setX}=\bar{\setK}\cap\setL$. 
    \item Unlike splitting methods or augmented Lagrangian algorithms \cite{birgin2017complexity,bolte2018nonconvex,grapiglia2020complexity,Andreani:2019uf,Andreani:2021vq,khanh2024new} our algorithms generate feasible approximate solutions and have complexity guarantees similar to the optimal complexity guarantees for unconstraied non-convex optimization.
    \item Since our algorithms generate a sequence of relatively interior points, the objective $f$ need not be differentiable at the relative boundary. A prominent example of such applications is the nonlinear regression problem with sparsity penalty:
    \begin{equation}
    \label{eq:l_p_regression}
        \min_{x \geq 0} \left\{f(x)=\ell(x)+\lambda \|x\|_p^p\right\},
    \end{equation}        
    where $\ell(x)$ is a non-convex loss function, $\lambda > 0$, $p \in (0,1)$.
    \item Our penalty-based approach provides a flexible framework since SCBs possess a calculus. Namely, a sum $h_{\bar{\setK}_1}+h_{\bar{\setK}_2}$ of SCBs for $\bar{\setK}_1$ and $\bar{\setK}_2$ is a SCB for $\bar{\setK}_1 \cap \bar{\setK}_2$. Moreover, SCBs can be efficiently constructed for a large variety of sets encountered in applications \cite{Nes18}. This is important, for example when $\bar{\setK}$ is given as an intersection of $1$-norm and Total Variation balls \cite{Liu2018ImageRU,hansen2023total}.
\end{enumerate}



\paragraph{Related works.}
Motivated, in particular, by training of neural networks, non-convex optimization is an active area of research in optimization and ML communities, see, e.g., the review \cite{danilova2022recent} and the references therein. An important part of this research concerns the global complexity guarantees of the proposed algorithms. 

\textbf{First-order methods.} 
If $f$ has Lipschitz gradient and there are no constraints, the standard gradient descent achieves the lower iteration complexity bound $O(\eps^{-2})$ to find a first-order $\eps$-stationary point $\hat{x}$ such that $\norm{\nabla f(\hat{x})} \leqslant\eps$  \cite{Nes18,CarDucHinSid19b,CarDucHinSid19}. 
In the composite optimization setting which includes problems with simple, projection-friendly, constraints a similar iteration complexity is achieved by the mirror descent algorithm \cite{lan2020first,ghadimi2016mini-batch,bogolubsky2016learning}. Various acceleration strategies of mirror and gradient descent methods have been derived in the literature, attaining the same bound as of gradient descent in the unconstrained case \cite{GhaLan16,guminov2019accelerated,nesterov2020primal-dual,guminov2021combination} or improving upon it under additional assumptions \cite{carmon2017convex,agarwal2017finding}. 
A potential drawback of these algorithms is the computationally expensive projection onto the set $\bar{\setX}=\bar{\setK}\cap \setL$.

\textbf{Second-order methods.}
If $f$ has Lipschitz Hessian and there are no constraints, cubic-regularized Newton method \cite{Gri81,NesPol06} and second-order trust region algorithms \cite{conn2000trust,cartis2012complexity,CurRobSam17} achieve the lower iteration complexity bound $O(\max\{\eps_1^{-3/2},\eps_2^{-3/2}\})$ \cite{CarDucHinSid19b,CarDucHinSid19} to find a second-order  $(\eps_1,\eps_2)$-stationary point $\hat{x}$ such that $\|\nabla f(\hat{x}) \|_2 \leq \eps_1$ and $\lambda_{\min} \left(\nabla^2 f(\hat{x})\right) \geq - \sqrt{\eps_2}$, where $\lambda_{\min}(\cdot)$ denotes the minimal eigenvalue of a matrix 
\footnote{A number of works, e.g. \cite{cartis2012complexity,NeiWr20}, consider an $(\eps_1,\eps_2)$-stationary point defined as $\hat{x}$ such that $\|\nabla f(\hat{x}) \|_2 \leq \eps_1$ and $\lambda_{\min} \left(\nabla^2 f(\hat{x})\right) \geq - \eps_2$ and the corresponding complexity $O(\max\{\eps_1^{-3/2},\eps_2^{-3}\})$. Our definition and complexity bound are the same up to the redefinition of $\eps_2$.}. 
Extensions for problems with simple projection-friendly constraints also exist \cite{CarGouToi12,birgin2017regularization,CarGouToi18} with the same iteration complexity bounds, as well as for problems with nonlinear equality and/or inequality constraints \cite{curtis2018complexity,hinder2018worst-case,cartis2019optimality,birgin2017complexity,grapiglia2020complexity,xie2019complexity}, but these works do not consider general set constraints as in \eqref{eq:Opt} and again may require projections on $\bar{\setX}=\bar{\setK}\cap \setL$.

\textbf{Barrier algorithms.} 
Existing barrier methods for non-convex optimization deal with some particular cases of \eqref{eq:Opt}, such as $\bar{\setK}$ being non-negative orthant \cite{Ye92,TseBomSch11,HBA-linear,BiaCheYe15,HaeLiuYe18,NeiWr20}, $\bar{\setK}$ being a symmetric cone \cite{He:2022aa,dvurechensky2021hessian}, $f$ being a quadratic function \cite{Ye92,FayLu06,LuYua07}. None of these works covers the general problem \eqref{eq:Opt}.

Summarizing, the existing works require one or two of the following crucial assumptions: a) Lipschitz continuity of the gradient and/or Hessian on the whole feasible set, b) no constraints, or simple convex constraints that allow an easy projection, or constraints that do not involve general feasible sets, e.g., only non-negativity constraints. Moreover, existing algorithms do not always come with complexity guarantees.
For further discussion, see Appendix \ref{sec:app_discussion}.

In this paper we develop a flexible and unifying algorithmic framework that is able to accommodate first- and second-order interior-point algorithms for \eqref{eq:Opt} with potentially non-convex and non-smooth at the relative boundary objective functions, and general set constraints. To the best of our knowledge, our framework is the first one providing complexity results for first- and second-order algorithms to reach points satisfying, respectively, suitably defined approximate first- and second-order necessary optimality conditions, under such weak assumptions and for such a general setting. 


\paragraph{Contributions.}
In this paper, we close the gap in the theoretical analysis of barrier algorithms for general non-convex objectives and general set constraints by constructing 
first- and second-order algorithms with complexities $O(\eps^{-1})$ and $O(\eps^{-3/2})$ respectively. In more detail, our contributions are as follows:

\textbf{Optimality conditions.} We propose a suitable set of first- and second-order necessary optimality conditions for \eqref{eq:Opt} that do not require $f$ to be differentiable at the relative boundary of $\bar{\setK}$. This is followed by the definition of approximate stationary points which we call $\eps$-KKT and $(\eps_1,\eps_2)$-2KKT points respectively (see Section \ref{sec:Optimality} for a precise definition).

\textbf{First-order algorithm.} 
We propose a new \emph{first-order adaptive barrier method} ($\AHBA$, Algorithm \ref{alg:AHBA}). 
To find a step direction, a linear model for $F_{\mu}$ in  \eqref{eq:potential}, regularized by the squared local norm induced by the Hessian of $h$, is minimized over the tangent space of the linear subspace $\setL$.
We prove that our new first-order method enjoys (see Theorem \ref{Th:AHBA_conv}) the upper iteration complexity bound $O(\eps^{-2})$ for reaching an $\eps$-KKT point when a ``descent Lemma'' holds relative to the local norm induced by the Hessian of $h$ (see Assumption \ref{ass:gradLip} for precise definition). Our algorithm is adaptive in the sense that it does not require the knowledge of the Lipschitz constant of the gradient.

\textbf{Second-order algorithm.} 
We propose a new \emph{second-order adaptive barrier method} ($\SAHBA$, Algorithm \ref{alg:SOAHBA}), for which the step direction is determined by a minimization subproblem over the same tangent space. 
But, in this case, the minimized model is composed of the linear model for $F_{\mu}$ augmented by second-order term for $f$ and regularized by the cube of the local norm induced by the Hessian of $h$. 
The regularization parameter is chosen adaptively in the spirit of \cite{NesPol06,CarGouToi12}.
The resulting minimization subproblem can be formulated as a non-convex optimization problem over a linear subspace, which can be solved as in the unconstrained case originally studied in \cite{NesPol06}. 
We establish (see Theorem \ref{Th:SOAHBA_conv}) the worst-case bound $O(\max\{\eps_1^{-3/2},\eps_2^{-3/2}\})$ on the number of iterations for reaching an $(\eps_1,\eps_2)$-2KKT point under a weaker version of assumption that the Hessian of $f$ is Lipschitz relative to the local norm induced by the Hessian of $h$ (see Assumption \ref{ass:2ndorder} for precise definition).

We do not perform numerical experiments for two reasons. 1) Our main goal is theoretical and we show that barrier algorithms possess favorable complexity beyond convexity. 2) We are not aware of any baseline algorithms with similar complexity that can accommodate with general set constraints, and produce feasible iterates without involving a projection step. 

\paragraph{Notation.} $\setE$ denotes a finite-dimensional real vector space, and $\setE^{\ast}$ the dual space, which is formed by all linear functions on $\setE$. The value of $s\in\setE^{\ast}$ at $x\in\setE$ is denoted by $\inner{s,x}$. The gradient and Hessian of a (twice) differentiable function $f:\setE\to\R$ at $x \in \setE$ are denoted as $\nabla f(x)\in\setE^{\ast}$ and $\nabla^{2}f(x)$ respectively. The directional derivative of function $f:\setE\to\R$ is defined in the usual way:
$Df(x)[v]=\lim_{\eps\to 0+} \frac{1}{\eps}[f(x+\eps v)-f(x)]$.
More generally, for $v_{1},\ldots,v_{p}\in\setE$, we define $D^{p}f(x)[v_{1},\ldots,v_{p}]$ the $p$-th directional derivative at $x$ along directions $v_{i}\in\setE$. In that way we define the gradient $\nabla f(x)\in\setE^{\ast}$ by $Df(x)[u]=\inner{\nabla f(x),u}$ and the Hessian $\nabla^{2}f(x):\setE\to\setE^{\ast}$ by $\inner{\nabla^{2}f(x)u,v}=D^{2}f(x)[u,v]$. 
For an operator $\bH:\setE\to\setE^{\ast}$, we denote by $\bH^{\ast}:\setE\to\setE^{\ast}$ its adjoint operator, defined by the identity 
$
(\forall u,v\in\setE): \qquad \inner{\bH u,v}=\inner{u,\bH^{\ast}v}.
$
For an operator $\bH:\setE\to\setE^{\ast}$ we say that it is positive semi-definite if $\inner{\bH u,u}\geq 0$ for all $ u\in\setE$, denoted as $\bH\succeq 0$. If the inequality is always strict for non-zero $u$, then $\bH$ is called positive definite denoted as $\bH\succ 0$. The standard Euclidean norm of a vector $x\in\setE$ is denoted by $\norm{x}$.
We denote by $\setL_{0}\eqdef\{ v \in \setE \vert \bA v = 0\}\eqdef\ker(\bA)$ the tangent space associated with the affine subspace $\setL\subset\setE$.

\paragraph{Self-concordant barriers.}
By our assumptions $\bar{\setK}\subset\setE$ has a self-concordant barrier $h(x)$ with finite parameter value $\nu$ \cite{NesNem94}.
\begin{definition}
\label{def:LHSCB}
A function $h:\bar{\setK}\to(-\infty,\infty]$ with $\dom h=\setK$ is called a $\nu$-\emph{self-concordant barrier} ($\nu$-SCB) for the set $\bar{\setK}$ if for all $x \in \setK$ and $u\in \setE$
\begin{align*}
&\abs{D^{3}h(x)[u,u,u]}\leq 2D^{2}h(x)[u,u]^{3/2},\text{ and } \\
&\sup_{u\in\setE}\abs{2 Dh(x)[u]-D^{2}h(x)[u,u]}\leq \nu.
\end{align*}
We denote the set of $\nu$-self-concordant barriers by $\scrH_{\nu}(\setK)$.
\end{definition}
Note that  $\nu\geq 1$ and the Hessian $H(x)\eqdef\nabla^{2}h(x):\setE\to\setE^{\ast}$ is a positive definite linear operator that defines the \emph{local norm} $\norm{u}_{x}\eqdef\inner{H(x)u,u}^{1/2}$.
The corresponding \emph{dual local norm} on $\setE^{\ast}$ is then defined as $\norm{s}_{x}^{\ast}\eqdef \inner{[H(x)]^{-1}s,s}^{1/2}$.

The \emph{Dikin ellipsoid}  with center $x\in\setK$ and radius $r>0$ is defined as 
$\scrW(x;r)\eqdef\{u\in\setE\vert\;\norm{u-x}_{x}<r\}$. This object allows us to guarantee the feasibility of the iterates in each iteration of our algorithms. 
\begin{lemma}[Theorem 5.1.5 \cite{Nes18}]
\label{lem:Dikin}
For all $x\in\setK$ we have $\scrW(x;1)\subseteq\setK$.
\end{lemma}
The following upper bound for the barrier $h$ is used to establish per-iteration decrease of the potential $F_{\mu}$.
\begin{proposition}[{Theorem 5.1.9 \cite{Nes18}}]
\label{prop:SCF_upper_bound}
Let $h\in\scrH_{\nu}(\setK)$, $x\in\dom h$, and a fixed direction $d \in \setE$. For all $t \in [0,\frac{1}{\norm{d}_x})$, with the convention that $\frac{1}{\norm{d}_{x}}=+\infty$ if $\norm{d}_x=0$, we have:
\begin{equation}
\label{eq:SCF_upper_bound} 
h(x + t d) \leq h(x) + t\inner{\nabla h(x),d} + t^2 \norm{d}_{x}^2 \omega(t \norm{d}_{x}),
\end{equation}
where $\omega(t)=\frac{-t-\ln(1-t)}{t^2}$.
\end{proposition}
\noindent
We will also use the following inequality for the function $\omega(t)$ \cite{Nes18}, Lemma 5.1.5:
\begin{equation}\label{eq:omega_upper_bound}
\omega(t) \leq \frac{1}{2(1-t)}, \; t \in [0,1).
\end{equation}
Appendix \ref{app:barrier} contains some more technical properties of SCBs which are used in the proofs.

\section{Approximate Optimality Conditions}
\label{sec:Optimality}
Given $\eps >0$, we define an $\eps$-approximate normal cone for the set $\bar{\setK}$ at $x \in \bar{\setK}$ \cite{burachik1997enlargement} to be the set 
\begin{equation}
\label{eq:approx_normal_cone_def}
\NC_{\bar{\setK}}^{\eps}(x) \eqdef \left\{ s \in \setE^* \vert \inner{s , y - x} \leq \eps \; \forall \; y \in \bar{\setK} \right\}.
\end{equation}
Clearly, we have $\NC_{\bar{\setK}}^{0}(x) = \NC_{\bar{\setK}}(x)$, where the latter denotes the normal cone for $\bar{\setK}$ at $x$.

The following result gives a necessary optimality condition for \eqref{eq:Opt}. Importantly, this result holds even if $x^*$ is at the boundary of $\bar{\setK}$ and $f$ is not differentiable at $x^*$.  
\begin{theorem}
\label{Th:limit_optimality}
Let the assumptions described above hold for problem \eqref{eq:Opt} and $x^* \in \bar{\setK}$ be a local solution to this problem.
Then, there exists a sequence of approximate solutions $x^k \in \setE$ and sequences of approximate Lagrange multipliers $y^k \in \R^m$,  $s^k \in \setE^*$ s.t.:
\begin{enumerate}
	\item $x^k \in \setK$, $\bA x^k= b$ for all $k$ and $x^k \to x^*$, \label{Th:limit_optimality_1}
	\item $\nabla f(x^k) -\bA^{\ast}y^{k} - s^k \to 0$, \label{Th:limit_optimality_2}
	\item $-s^k \in \NC_{\bar{\setK}}^{\sigma_k}(x^{k})$, where $\sigma_k \to 0$. \label{Th:limit_optimality_3} 
\end{enumerate}
If in addition $f$ is twice differentiable on $ \setK$, then there exist $\theta_{k},\delta_{k}\in(0,\infty)$ such that $\theta_{k},\delta_{k}\to 0$ and 
\begin{equation}\label{Th:limit_optimality_4}
 \inner{(\nabla^2 f(x^k) + \theta_k H(x^k) +\delta_{k}\bI)d,d}\geq 0
\end{equation}
for all $d \in \setE$ such that $\bA d=0$. Here $\bI$ -- identity operator.
\end{theorem}

The proof is based on interior-penalty arguments sketched below, and fully detailed in Appendix \ref{sec:optimality_proof}. Namely, there exists a sequence $x^k \to x^*$ that solves the sequence of penalized problems 
\begin{equation}
\label{eq:main_Haeser_limit_optimality_proof_3}
	\min_{x}  f(x) + \frac{1}{4}\|x-x^*\|^4 + \mu_k h(x)  \quad  \text{s.t.: } \; \bA x = b,
\end{equation}
where $\mu_k> 0$, $\mu_k \to 0$. The first- and second-order optimality conditions for the latter problem imply then the statement of the Theorem. 

The most interesting part of the above result for us is the second-order condition \eqref{Th:limit_optimality_4} since it allows us to certify second-order stationary points using the Hessian $H(x)$ of the barrier. At the same time, the first-order conditions may be strengthened compared to the ones in Theorem \ref{Th:limit_optimality}.
Indeed, if $x^{\ast}$ is a local solution of the optimization problem \eqref{eq:Opt} at which the objective function $f$ is continuously differentiable, then there exists $y^{\ast}\in\R^{m}$ such that $
\nabla f(x^{\ast})-\bA^{\ast}y^{\ast}\in-\NC_{\bar{\setK}}(x^{\ast})$, 
or, equivalently,
\begin{equation}\label{eq:Fermat_2}
\inner{\nabla f(x^{\ast})-\bA^{\ast}y^{\ast},x-x^{\ast}} \geq 0 \qquad\forall x\in\bar{\setK}. 
\end{equation}
The standard way to construct an approximate first-order optimality condition is to add an $\eps$-perturbation in the r.h.s. of \eqref{eq:Fermat_2}:
\begin{equation}\label{eq:Appr_Fermat_2}
\inner{\nabla f(\bar{x})-\bA^{\ast}\bar{y},x-\bar{x}} \geq -\eps \qquad\forall x\in\bar{\setK}. 
\end{equation}
which is equivalent to  $-(\nabla f(\bar{x})-\bA^{\ast}\bar{y})\in\NC_{\bar{\setK}}^{\eps}(\bar{x})$. 

Motivated by the above, we introduce the following notion of an approximate first-order KKT point for problem \eqref{eq:Opt}.
\begin{definition}\label{def:eps_KKT}
Given $\eps \geq 0$, a point $\bar{x}\in\setE$ is an $\eps$-KKT point for problem \eqref{eq:Opt} if there exists $\bar{y}\in\R^{m}$ such that
\begin{align}
&\bA\bar{x}=b,\bar{x}\in\setK,\label{eq:eps_optim_equality_set}\\
&\inner{\nabla f(\bar{x})-\bA^{\ast}\bar{y},x-\bar{x}} \geq -\eps \qquad\forall x\in\bar{\setK}.\label{eq:eps_optim_grad}
\end{align}
\end{definition}

We underline that when $\eps_k \to 0$, every convergent subsequence $(x^{k_j})_{j\geq 1}$ of a sequence $(x^{k})_{k\geq 1}$ of $\eps_k$-KKT points converges to a stationary point in the sense of Theorem \ref{Th:limit_optimality}.
Indeed, clearly, such subsequence satisfies item  \ref{Th:limit_optimality_1}. After defining $s^{k_j}=\nabla f(x^{k_j})-\bA^{\ast} y^{k_j}$, we see that item \ref{Th:limit_optimality_2} trivially holds as equality $\nabla f(x^{k_j})-\bA^{\ast} y^{k_j}-s^{k_j}=0$. Finally, the definition of $s^{k_j}$, \eqref{eq:eps_optim_grad}, \eqref{eq:approx_normal_cone_def}, and the condition $\eps_k \to 0$ imply item \ref{Th:limit_optimality_3} with $\sigma_{k_j}=\eps_{k_j}$. Thus, the limit of the subsequence $x^{k_j}$ satisfies the first three items of Theorem \ref{Th:limit_optimality}, and thus is a first-order stationary point according to this theorem.

Based on the second-order condition \eqref{Th:limit_optimality_4} in Theorem \ref{Th:limit_optimality}, we can augment Definition \ref{def:eps_KKT} with an approximate second-order condition. This leads us to the following notion of an approximate second-order KKT point for problem \eqref{eq:Opt}.
\begin{definition}\label{def:eps_SOKKT}
Given $\eps_1,\eps_2 \geq 0$, a point $\bar{x}\in \setE$ is an $(\eps_1,\eps_2)$-2KKT point for problem \eqref{eq:Opt} if there exists $\bar{y}\in\R^{m}$ such that 
\begin{align}
& \bA\bar{x}=b, \bar{x}\in \setK, \label{eq:eps_SO_optim_equality_set} \\
& \inner{\nabla f(\bar{x})-\bA^{\ast}\bar{y},x-\bar{x}} \geq -\eps_1 \qquad\forall x\in\bar{\setK},  \label{eq:eps_SO_optim_FO} \\
& \nabla^2f(\bar{x}) +\sqrt{\eps_{2}} H(\bar{x}) \succeq 0 \;\; \text{on} \;\; \setL_0\rev{= \{ v \in \setE \vert \bA v = 0\}}.   \label{eq:eps_SO_optim_SO}
\end{align}
\end{definition}
Note that our definition of an approximate second-order KKT point is motivated by the notion of weak second-order approximate stationary conditions for non-convex optimization using barrier algorithms \cite{HaeLiuYe18,NeiWr20,He:2022aa}. Just as for Definition \ref{def:eps_KKT}, we can prove that every accumulation point of a sequence of  $(\eps_{1,k},\eps_{2,k})$-2KKT points satisfies items \ref{Th:limit_optimality_1}, \ref{Th:limit_optimality_2}, \ref{Th:limit_optimality_3} of Theorem \ref{Th:limit_optimality}. Setting $\theta_{k_j}= \sqrt{\eps_{2,k_j}}$ and $\delta_k=0$, we see that the condition \eqref{Th:limit_optimality_4} also holds. Thus, the limit of the subsequence $x^{k_j}$ satisfies all the four items of Theorem \ref{Th:limit_optimality}, and thus is a second-order stationary point according to this theorem. An important advantage of the above definitions is that $\bar{x}$ lies in the relative interior of the feasible set. Thus, $f$ may be non-differentiable at the relative boundary of the feasible set, see, e.g., problem \eqref{eq:l_p_regression}. Moreover, we can use the Hessian of the barrier $H(\bar{x})$ in the second-order condition since $\bar{x}$ is in the interior of $\bar{\setK}$.

\section{First-Order Barrier Algorithm}
\label{sec:firstorder}

In this section we introduce our first-order potential reduction method for solving \eqref{eq:Opt} that uses a barrier $h \in\scrH_{\nu}(\setK)$ and potential function \eqref{eq:potential}. 

\subsection{Smoothness Assumption}
Given $x\in\setX$, define the set of \emph{feasible directions} as $\scrF_{x}\eqdef \{v\in\setE\vert x+v\in\feas\}.$ Lemma \ref{lem:Dikin} implies that 
\begin{equation}\label{eq:Dikinv}
\scrT_{x}\eqdef\{v\in\setE\vert \bA v=0,\norm{v}_{x}<1\}\subseteq\scrF_{x}.
\end{equation}
Upon defining $d=[H(x)]^{1/2}v$ for $v\in\scrT_{x}$, we obtain a direction $d$ satisfying $\bA[H(x)]^{-1/2}d=0$ and $\norm{d}=\norm{v}_{x}$. Hence, for $x\in\setK$, we can equivalently characterize the set $\scrT_{x}$ as $\scrT_{x}=\{[H(x)]^{-1/2}d\vert \bA[H(x)]^{-1/2}d=0,\norm{d}<1\}$.\\ 
For the analysis of the first-order algorithm we use the following first-order smoothness condition.
\begin{assumption}[Local smoothness]
\label{ass:gradLip}
$f:\setE\to\R\cup\{+\infty\}$ is continuously differentiable on $\feas$ and there exists a constant $M>0$ such that for all $x\in\feas$ and $v\in\scrT_{x}$, where $\scrT_{x}$ is defined in \eqref{eq:Dikinv}, we have 
\begin{equation}\label{eq:gradLip}
f(x+v) - f(x) - \inner{\nabla f(x),v} \leq \frac{M}{2}\norm{v}_x^2.
\end{equation}
\end{assumption}
\begin{remark}
\label{rem:bounded1}
If the set $\bar{\setX}$ is bounded, we have $\lambda_{\min}(H(x)) \geq \sigma$ for some $\sigma >0$. In this case, assuming $f$ has an $M$-Lipschitz continuous gradient, the classical descent lemma \cite{Nes18} implies Assumption \ref{ass:gradLip}. Indeed,
\[
f(x+v) - f(x) - \inner{\nabla f(x),v} \leq \frac{M}{2}\norm{v}^2 \leq \frac{M}{2\sigma}\norm{v}_x^2. \;\;\;\text{\close}
\]
\end{remark}
Considering $x\in\setX,v\in\scrT_{x}$ and combining eq. \eqref{eq:gradLip} with eq. \eqref{eq:SCF_upper_bound} (with $d=v$ and $t=1< \frac{1}{\norm{v}_x}$) gives us the following upper bound that holds for all $x\in\setX,v\in\scrT_{x}$ and $L\geq M$  
\begin{align}\label{eq:descentFOM}
F_{\mu}(x+v)\leq &F_{\mu}(x)+\inner{\nabla F_{\mu}(x),v}+\frac{L}{2}\norm{v}^{2}_{x}\nonumber \\
&+\mu\norm{v}^{2}_{x}\omega(\norm{v}_x).
\end{align}

\subsection{Algorithm and Its Complexity}
\label{S:FO_descr}
We assume that our algorithm starts from a $\nu$-analytic center, i.e. a point $x^{0}\in\setX$ such that 
\begin{equation}\label{eq:analytic_center}
h(x)\geq h(x^{0})-\nu\qquad\forall x\in\setX. 
\end{equation}
Obtaining such a point requires solving a \textit{convex} optimization problem $\min_{x\in \feas}h(x)$ up to a very loose accuracy $\nu \geq 1$.  
We denote $\Delta^f_0 \eqdef f(x^{0}) - f_{\min}(\setX)$.

\paragraph{Defining the  search direction.}
Let $x \in \feas$ be given. Our first-order method uses a quadratic model
\[
Q^{(1)}_{\mu}(x,v) \eqdef F_{\mu}(x) + \inner{\nabla F_{\mu}(x),v}+\frac{1}{2}\norm{v}_{x}^{2} 
\]
to compute a search direction $v_{\mu}(x)$, given by 
\begin{equation}\label{eq:search}
v_{\mu}(x) \eqdef \argmin_{v\in\setE:\bA v=0}    Q^{(1)}_{\mu}(x,v) .
\end{equation}
This search direction is determined by the following system of optimality conditions involving the dual variable $y_{\mu}(x)\in\R^{m}$:
\begin{align}
\nabla F_{\mu}(x) + H(x)v_{\mu}(x) - \bA^{\ast} y_{\mu}(x) &= 0, \label{eq:finder_1} \\
\bA  v_{\mu}(x) &=0. \label{eq:finder_2}
\end{align}
Since $H(x)\succ 0$ for $x\in\feas$, any standard solution method \citep{NocWri00} can be applied for the above linear system.
Since $H(x)\succ 0$ for $x\in\feas$, and $\bA$ has full row rank, the optimality conditions have a unique solution.

\paragraph{Defining the step-size.}
Consider a point $x \in\setX$ and a point $x^{+}(t)\eqdef x+tv_{\mu}(x)$, where $t\geq 0$ is the step-size. Our aim is to choose $t$ to ensure the feasibility of iterates and decrease of the potential. By Lemma \ref{lem:Dikin} and \eqref{eq:finder_2}, we know that $x^{+}(t)\in\feas$ for all $t\in I_{x,\mu} \eqdef [0,\frac{1}{\norm{v_{\mu}(x)}_{x}})$. Multiplying \eqref{eq:finder_1} by $v_{\mu}(x)$ and using \eqref{eq:finder_2}, we obtain 
$\inner{\nabla F_{\mu}(x),v_{\mu}(x)}=-\norm{v_{\mu}(x)}_{x}^{2}$. Choosing $t \in  I_{x,\mu}$, we have
\[
t^{2}\norm{v_{\mu}(x)}_{x}^{2}\omega(t\norm{v_{\mu}(x)}_{x}) \stackrel{\eqref{eq:omega_upper_bound}}{\leq}  \frac{t^{2}\norm{v_{\mu}(x)}_{x}^{2}}{2(1-t\norm{v_{\mu}(x)}_{x})}. 
\]
Therefore, if $t\norm{v_{\mu}(x)}_{x}\leq 1/2$, we get from \eqref{eq:descentFOM} that
\begin{align}
&F_{\mu}(x^{+}(t))-F_{\mu}(x) \nonumber \\
&\leq -t\norm{v_{\mu}(x)}_{x}^{2}+\frac{t^{2}M}{2}\norm{v_{\mu}(x)}_{x}^{2}+\mu t^{2}\norm{v_{\mu}(x)}_{x}^{2} \nonumber\\
&= -t \norm{v_{\mu}(x)}_{x}^{2}\left(1-\frac{M+2\mu}{2}t\right) \eqdef -\eta_{x}(t).\label{eq:success}
\end{align}
The function $\eta_{x}(t)$ is strictly concave with the unique maximum at $ \frac{1}{M+2\mu}$. 
Thus, maximizing the per-iteration decrease $\eta_{x}(t)$  under the restriction $0\leq t\leq\frac{1}{2\norm{v_{\mu}(x)}_{x}}$, we choose the step-size
\begin{equation*}
\ct_{\mu,M}(x)\eqdef \min \left\{\frac{1}{M+2\mu},\frac{1}{2\norm{v_{\mu}(x)}_{x}}\right\}.
\end{equation*}
\paragraph{Adaptivity to the Lipschitz constant.}
To get rid of the explicit dependence of the step size on the Lipschitz parameter $M$, we propose a backtracking/adaptive procedure in the spirit of \cite{NesPol06}. This procedure generates a sequence of positive numbers $(L_{k})_{k\geq 0}$ for which the local Lipschitz smoothness condition \eqref{eq:gradLip} holds. More specifically, let $x^{k}$ be the current position of the algorithm with the corresponding initial local Lipschitz estimate $L_{k}$ and $v^{k}=v_{\mu}(x^{k})$ is the corresponding search direction. To determine the next iterate $x^{k+1}$, we iteratively try step-sizes $\alpha_k$ of the form $\ct_{\mu,2^{i_k}L_k}(x^{k})$ for $i_k\geq 0$ until the local smoothness condition \eqref{eq:gradLip} holds with $x=x^{k}$, $v= \alpha_k v^{k}$ and local Lipschitz estimate $M=2^{i_k}L_k$, see \eqref{eq:LS}. This process must terminate in finitely many steps 
 since when $2^{i_k}L_k \geq M$, inequality \eqref{eq:gradLip} with $M$ changed to $2^{i_k}L_k$, i.e., \eqref{eq:LS}, follows from Assumption \ref{ass:gradLip}. 
Setting $L_{k+1} = 2^{i_k-1}L_{k}$ allows $L_k$ to adaptively decrease so that in the areas where $f$ is more smooth the algorithm uses larger step-sizes.

Combining the definition of the search direction  \eqref{eq:search} with the above backtracking strategy, yields a  First-Order Adaptive Barrier Method ($\AHBA$, Algorithm \ref{alg:AHBA}).
\begin{algorithm}[t]
\caption{First-Order Adaptive Barrier Method  - $\AHBA(\mu,\eps,L_{0},x^{0})$}
\label{alg:AHBA}
\SetAlgoLined
\KwData{ $h \in\scrH_{\nu}(\setK)$, $\mu>0,\eps>0,L_0>0,x^{0}\in\setX$.
}
\KwResult{$(x^{k},y^{k},s^{k},L_{k})\in\setX\times\R^{m}\times\setE^{\ast}\times\R_{+}$, where $s^{k}=\nabla f(x^{k}) -\bA^{\ast}y^{k}$, and $L_{k}$ is the last estimate of the Lipschitz constant.}
Set $k=0$\;
\Repeat{
		$\norm{v^k}_{x^k} < \tfrac{\eps}{3\nu}$ 
	}{	
		Set $i_k=0$. Find $v^k\eqdef v_{\mu}(x^k)$ and  the corresponding dual variable $y^k\eqdef y_{\mu}(x^k)$ as the solution to
		\begin{equation}\label{eq:finder}
		\min_{v\in\setE:\bA v=0}\{F_{\mu}(x^k) + \inner{\nabla F_{\mu}(x^k),v}+\frac{1}{2}\norm{v}_{x^{k}}^{2}\}. 
		\end{equation}
		\Repeat{
			\begin{equation}				
				\hspace{-3em}f(z^{k}) \leq f(x^{k}) + \inner{\nabla f(x^{k}),z^{k}-x^{k}}+2^{i_{k}-1}L_{k}\norm{z^{k}-x^{k}}^{2}_{x^{k}}.
				\label{eq:LS}
			\end{equation}
		}
		{
			\begin{equation}\label{eq:alpha_k}
				\hspace{-1.5em}\alpha_k \eqdef  \min \left\{\frac{1}{2^{i_k}L_{k} + 2 \mu},\frac{1}{2\norm{v^k}_{x^{k}}} \right\}	
			\end{equation}
			Set $z^{k}=x^{k} + \alpha_k v^k$, $i_k=i_k+1$\;
		}
		Set $L_{k+1} = 2^{i_k-1}L_{k}$, $x^{k+1}=z^{k}$, $k=k+1$\;
	}
\end{algorithm}

\paragraph{Complexity bound.}
Our main result on the iteration complexity of Algorithm \ref{alg:AHBA} is the following Theorem. 
\begin{theorem}
\label{Th:AHBA_conv}
Let Assumptions \ref{ass:1} and \ref{ass:gradLip} hold. Fix the error tolerance $\eps>0$, the regularization parameter $\mu=\frac{\eps}{\nu}$, and some initial guess $L_0>0$ for the Lipschitz constant in \eqref{eq:gradLip}. Let $(x^{k})_{k\geq 0}$ be the trajectory generated by $\AHBA(\mu,\eps,L_{0},x^{0})$, where $x^{0}$ is a $\nu$-analytic center satisfying \eqref{eq:analytic_center}. Then the algorithm stops in no more than 
\begin{equation}
\label{eq:FO_main_Th_compl}
\hspace{-1em}\K_{I}(\eps,x^{0})= \ceil[\bigg]{36(\Delta^f_0+ \eps) \frac{\nu^{2}(\max\{M,L_0\}+\eps/\nu)}{\eps^{2}}}
\end{equation}
outer iterations, and the number of inner iterations is no more than $2(\K_{I}(\eps,x^{0})+1)+\max\{\log_{2}(M/L_{0}),0\}$. Moreover, the last iterate obtained from $\AHBA(\mu,\eps,L_{0},x^{0})$ constitutes a $2\eps$-KKT point for problem \eqref{eq:Opt} in the sense of Definition \ref{def:eps_KKT}.
\end{theorem}
\begin{remark}
\label{rem:FO_complexity_simplified}
Since $\nu \geq 1$, $\Delta^f_0$ is expected to be larger than $\eps$, and the constant $M$ is potentially large, we see that the main term in the complexity bound \eqref{eq:FO_main_Th_compl} is $O\left(\frac{M\nu^2\Delta^f_0}{\eps^2}\right)=O\left(\frac{1}{\eps^{2}}\right)$, i.e., has the same dependence on $\eps$ as the standard complexity bounds  \cite{CarDucHinSid19b,CarDucHinSid19,lan2020first} of first-order methods for non-convex problems under the standard Lipschitz-gradient assumption, which on bounded sets is subsumed by our Assumption \ref{ass:gradLip}. Further, if the function $f$ is linear, Assumption \ref{ass:gradLip} holds with $M=0$ and we can take $L_0=0$. In this case, the complexity bound \eqref{eq:FO_main_Th_compl} improves to $O\left(\frac{\nu\Delta^f_0}{\eps}\right)$. \close
\end{remark}
\begin{remark}
\label{Rm:FO_restarts}
A potential drawback of Algorithm \ref{alg:AHBA} may be that it requires to fix the parameter $\eps$ before start. This may be easily resolved by a restart procedure with warm starts. Namely, we run Algorithm \ref{alg:AHBA} in epochs numbered by $i\geq 0$. For each restart, we choose $x^0$ as the output of the previous restart and set $\eps=\eps_i=2^{-i}\eps_0$. Such an algorithm may be run infinitely. To reach any desired accuray $\eps$, it is sufficient to make $I=I(\eps)=O(\log_{2}(\eps_{0}/\eps))$ restarts. Then, the total number of inner iterations is of the order $\sum_{i=0}^{I-1} \frac{C}{\eps_i^2}=\sum_{i=0}^{I-1} \frac{C2^{i}}{\eps_0^2} = O(\eps^{-2})$, i.e., the complexity is the same. \close
\end{remark}
\paragraph{Sketch of the proof of Theorem \ref{Th:AHBA_conv}.}
The proof is organized in three steps. First, we prove the correctness of the algorithm, i.e., that it generates a sequence of points in $\setX$, and, thus, is indeed an interior-point method. This follows from the construction of the algorithm by an induction argument. Next, we show that the line-search process of finding appropriate $L_k$ in each iteration is finite, and estimate the total number of trials in this process. 
Then we enter the core of our analysis where we prove that, if the stopping criterion does not hold at iteration $k$, i.e., $\norm{v^{k}}_{x^{k}} \geq \tfrac{\eps}{3\nu}$, then the objective $f$ is decreased by a quantity $O(\eps^{2})$, which follows from \eqref{eq:success}. Since the objective is globally lower bounded, we conclude that the method stops in at most $O(\eps^{-2})$ iterations. Finally, we show that when the stopping criterion holds, i.e., $\norm{v^{k}}_{x^{k}} < \tfrac{\eps}{3\nu}$, the method has generated an $\eps$-KKT point. On a high level, the result follows from the optimality condition \eqref{eq:finder_1} which implies by the stopping condition 
\[
\nabla f(x^{k}) - \bA^{\ast} y^{k} =  - H(x^{k})v^{k} - \mu \nabla h(x^{k}) = O(\eps + \mu \nu).
\]
The full proof is deferred to Appendix \ref{sec:App_FO}.

\section{Second-Order Barrier Algorithm}
\label{sec:secondorder}
In this section, we present a second-order method that uses also the Hessian of $f$ and the following assumption. 
\begin{assumption}[Local second-order smoothness]
\label{ass:2ndorder}
$f:\setE\to\R\cup\{+\infty\}$ is twice continuously differentiable on $\feas$ and there exists a constant $M>0$ such that, for all $x\in\feas$ and $v\in\scrT_{x}$, where $\scrT_{x}$ is defined in \eqref{eq:Dikinv}, we have
\begin{equation}\label{eq:SO_Lipschitz_Gradient}
\norm{\nabla f(x+v)-\nabla f(x)-\nabla^{2}f(x)v}^{\ast}_{x}\leq\frac{M}{2}\norm{v}^{2}_{x}.
\end{equation}
\end{assumption}
By integration, it is easy to show that a sufficient condition for \eqref{eq:SO_Lipschitz_Gradient} is the following local counterpart of the global Lipschitz condition on
$\nabla^2f$ \cite{Nes18}:
\begin{equation}\label{eq:LipHess}
\norm{\nabla^{2}f(x+u)-\nabla^{2}f(x+v)}_{\text{op},x}\leq M\norm{u-v}_{x},
\end{equation}
where 
$\norm{\BB}_{\text{op},x}\eqdef\sup_{u:\norm{u}_{x}\leq 1}\left\{\frac{\norm{\BB u}_{x}^{\ast}}{\norm{u}_{x}}\right\}$ is the induced operator norm for a linear operator $\BB:\setE\to\setE^{\ast}$. 
Further, again by integration \eqref{eq:SO_Lipschitz_Gradient} implies
\begin{align}\label{eq:cubicestimate}
&\hspace{-1em}f(x+v)-\left[f(x)+\inner{\nabla f(x),v}+\frac{1}{2}\inner{\nabla^{2}f(x)v,v}\right] \nonumber \\
& \leq\frac{M}{6}\norm{v}^{3}_{x}.
\end{align}
\begin{remark}
\label{Rm:Hess_Lip}
Assumption \ref{ass:2ndorder} subsumes, if $\bar{\setX}$ is bounded, the standard Lipschitz-Hessian setting. If the Hessian of $f$ is $M$-Lipschitz w.r.t. the standard Euclidean norm, we have by \cite{Nes18}, Lemma 1.2.4, that 
\[
\norm{\nabla f(x+v)-\nabla f(x)-\nabla^{2}f(x)v}\leq\frac{M}{2}\norm{v}^{2}.
\]
Since $\bar{\setX}$ is bounded, one can observe that $\lambda_{\max}([H(x)]^{-1})^{-1}=\lambda_{\min}(H(x)) \geq \sigma$ for some $\sigma >0$, and \eqref{eq:SO_Lipschitz_Gradient} holds. Indeed, denoting $g=\nabla f(x+v)-\nabla f(x)-\nabla^{2}f(x)v$, we obtain
\begin{align*}
    (\norm{g}_x^*)^2 &\leq \lambda_{\max}([H(x)]^{-1})\norm{g}^{2} \leq \frac{M^2}{4\lambda_{\min}(H(x))}\norm{v}^{4} \\
    &\leq \frac{M^2}{4\sigma^{3}}\norm{v}_x^4. \qquad\qquad\qquad\qquad\qquad\qquad\qquad\text{\close}
\end{align*}
\end{remark}
Assumption \ref{ass:2ndorder} also implies, via \eqref{eq:cubicestimate} and \eqref{eq:SCF_upper_bound} (with $d=v$ and $t=1< \frac{1}{\norm{v}_x}$), the following upper bound for $F_{\mu}$ that holds for all $x\in\setX,v\in\scrT_{x}$ and $L\geq M$: 
\begin{align}\label{eq:cubicdecrease}
F_{\mu}(x+v)\leq &F_{\mu}(x)+\inner{\nabla F_{\mu}(x),v}+\frac{1}{2}\inner{\nabla^{2}f(x)v,v}\nonumber\\
&+\frac{L}{6}\norm{v}^{3}_{x}+\mu\norm{v}^{2}_{x}\omega(\norm{v}_{x}).
\end{align}

\subsection{Algorithm and Its Complexity}
\label{S:SO_alg_descr}
\paragraph{Defining the search direction.} Let $x\in\setX$ be given. In order to find a search direction, we choose a parameter $L>0$, construct a cubic-regularized model of the potential $F_{\mu}$ \eqref{eq:potential}
\begin{align}\label{eq:Q_2}
Q^{(2)}_{\mu,L}(x,v)\eqdef &F_{\mu}(x)+\inner{\nabla F_{\mu}(x),v}+\frac{1}{2}\inner{\nabla^{2}f(x)v,v} \nonumber \\
&+\frac{L}{6}\norm{v}_{x}^{3}, 
\end{align}
and minimize it on the linear subspace $\setL_{0}$:
\begin{align}\label{eq:cubicproblem}
v_{\mu,L}(x)\in\Argmin_{v\in\setE:\bA v=0}  Q^{(2)}_{\mu,L}(x,v),
\end{align}
where by $\Argmin$ we denote the set of global minimizers. The model consists of three parts: linear approximation of $h$, quadratic approximation of $f$, and a cubic regularizer with penalty parameter $L>0$. Since this model and our algorithm use the second derivative of $f$, we call it a second-order method.
Our further derivations rely on the first-order optimality conditions for the problem \eqref{eq:cubicproblem}, which say that there exists $y_{\mu,L}(x)\in\R^{m}$ such that $v_{\mu,L}(x)$ satisfies
\begin{align}
\hspace{-5em}\nabla F_{\mu}(x)+\nabla^{2}f(x)v_{\mu,L}(x)& \nonumber \\
+\frac{L}{2}\norm{v_{\mu,L}(x)}_{x}H(x)v_{\mu,L}(x) - \bA^{\ast} y_{\mu,L}(x) &= 0,\label{eq:opt1}\\
 - \bA v_{\mu,L}(x)&=0. \label{eq:opt2}
\end{align}
We also use the following extension of \cite{NesPol06}, Prop. 1, with the local norm induced by $H(x)$.
\begin{proposition}
\label{Th:PD}
For all $x\in\feas$ it holds 
\begin{equation}\label{eq:PD}
\nabla^{2}f(x)+\frac{L}{2}\norm{v_{\mu,L}(x)}_{x}H(x)\succeq 0\qquad\text{ on }\;\;\setL_{0}.
\end{equation}
\end{proposition}

\paragraph{Defining the step-size.} 
To define the step-size, we act in the same fashion as in Section \ref{sec:firstorder} by considering $x^{+}(t)\eqdef x+t v_{\mu,L}(x)$, where $t\geq 0$ is a step-size. Using the optimality conditions \eqref{eq:opt1}, \eqref{eq:opt2}, \eqref{eq:PD}, we estimate (the full derivation is given in Appendix \ref{sec:App_SO}) the progress parameterized by $t$:
\begin{align}
&F_{\mu}(x^{+}(t))-F_{\mu}(x) \nonumber\\
&\leq - \frac{Lt^2\norm{v_{\mu,L}(x)}^{3}_{x} }{12} \left(3 - 2t \right)+\mu t^{2}\norm{v_{\mu,L}(x)}_{x}^{2} \nonumber \\
&= - \norm{v_{\mu,L}(x)}^{3}_{x} \frac{Lt^2}{12}\left(3 - 2t - \frac{12 \mu}{L\norm{v_{\mu,L}(x)}_{x}} \right) \eqdef -\eta_{x}(t).\label{eq:progressSOM_main}
\end{align}
The above inequality holds for all  $t \geq 0$ s.t. $t\norm{v_{\mu,L}(x)}_{x}\leq 1/2$.
To respect these constraints and guarantee that the decrease is positive, we choose the following step-size rule
\begin{equation}
\label{eq:SO_t_opt_def_main}
\ct_{\mu,L}(x) \eqdef  \min \left\{1, \frac{1}{2\norm{v_{\mu,L}(x)}_{x}}  \right\}.
\end{equation}  
Note that $\ct_{\mu,L}(x) \leq 1$ and $\ct_{\mu,L}(x)\norm{v_{\mu,L}(x)}_{x}\leq 1/2$. Thus, this choice of the step-size is feasible to derive \eqref{eq:progressSOM_main}. 

\paragraph{Adaptivity to the Lipschitz constant.} Just like Algorithm \ref{alg:AHBA}, our second-order method employs a line-search procedure to estimate the Lipschitz constant $M$ in \eqref{eq:SO_Lipschitz_Gradient}, \eqref{eq:cubicestimate} in the spirit of \cite{NesPol06,CarGouToi12}. More specifically, suppose that $x^{k}\in\setX$ is the current position of the algorithm with the corresponding initial local Lipschitz estimate $M_{k}$. To determine the next iterate $x^{k+1}$, we solve problem \eqref{eq:cubicproblem} with $L= L_k = 2^{i_k}M_{k}$ starting with $ i_k =0$, find the corresponding search direction $v^{k}=v_{\mu,L_{k}}(x^{k})$ and the new point $x^{k+1} = x^{k} + \ct_{\mu, L_{k}}(x^{k})v^{k}$. Then, we check whether the inequalities \eqref{eq:SO_Lipschitz_Gradient} and \eqref{eq:cubicestimate} hold with  $M=L_{k}$, $x=x^{k}$, $v = \ct_{\mu, L_{k}}(x^{k})v^{k}$, see \eqref{eq:SO_LS_2} and \eqref{eq:SO_LS_1}. If they hold, we make a step to $x^{k+1}$. Otherwise, we increase $i_k$ by 1 and repeat the procedure. Obviously, when $L_{k} = 2^{i_k}M_k \geq M$, both inequalities \eqref{eq:SO_Lipschitz_Gradient} and \eqref{eq:cubicestimate} with $M$ changed to $L_k$, i.e., \eqref{eq:SO_LS_2} and \eqref{eq:SO_LS_1}, are satisfied and the line-search procedure ends. For the next iteration we set $M_{k+1} = \max\{2^{i_k-1}M_{k},\underline{L}\}=\max\{L_{k}/2,\underline{L}\}$, so that the estimate for the local Lipschitz constant on the one hand can decrease allowing larger step-sizes, and on the other hand is bounded from below. 

\paragraph{Algorithm.} Combining the definition of the search direction in \eqref{eq:cubicproblem} with the just outlined backtracking strategy, yields a  Second-Order Adaptive Barrier Method ($\SAHBA$, Algorithm \ref{alg:SOAHBA}).
\begin{algorithm}[t!]
\caption{Second-Order Adaptive Barrier Method - $\SAHBA(\mu,\eps,M_{0},x^{0})$}
\label{alg:SOAHBA}
\SetAlgoLined
\KwData{ $h \in\scrH_{\nu}(\setK)$, $\mu>0,\eps>0,M_0\geq 144  \eps,x^{0}\in\setX$.}
\KwResult{$(x^{k},y^{k-1},s^{k},M_{k})\in\setX\times\R^{m}\times\setK^{\ast}\times\R_{+}$, where $s^{k}=\nabla f(x^{k}) -\bA^{\ast}y^{k-1}$, and $M_{k}$ is the last estimate of the Lipschitz constant.}
Set $\underline{L} \eqdef 144  \eps$, $k=0$\;
\Repeat{
		$\norm{v^{k-1}}_{x^{k-1}} < \Delta_{k-1}\eqdef\sqrt{\frac{\eps}{12L_{k-1}\nu}}$ \textnormal{ and }$\|v^{k}\|_{x^{k}} < \Delta_{k}\eqdef\sqrt{\frac{\eps}{12L_{k}\nu}}$
	}{	
		Set $i_k=0$.
		
		\Repeat{
		\begin{align}
			& \hspace{-3em} f(z^{k}) \leq f(x^{k}) + \inner{\nabla f(x^{k}),z^{k}-x^{k}} \nonumber\\
        &\hspace{-2em}+\frac{1}{2}\inner{\nabla^{2}f(x^{k})(z^{k}-x^{k}),z^{k}-x^{k}} +\frac{L_{k}}{6}\norm{z^{k}-x^{k}}^{3}_{x^{k}}, \label{eq:SO_LS_1}\\
		& \textnormal{and} \;\; \norm{\nabla f(z^{k})-\nabla f(x^{k})-\nabla^{2}f(x^{k})(z^{k}-x^{k})}^{\ast}_{x^{k}} \nonumber\\
        &\leq\frac{L_{k}}{2}\norm{z^{k}-x^{k}}^{2}_{x^{k}}. \label{eq:SO_LS_2}
		\end{align}
		}
		{
			Set $L_k = 2^{i_k}M_k$. Find $v^k \eqdef v_{\mu,L_k}(x^{k})$ and $y^k \eqdef y_{\mu,L_k}(x^{k})$ as a global solution to
			\begin{align}
                &  \hspace{-2em} \min_{v:\bA v=0} Q^{(2)}_{\mu,L_k}(x^k,v), \;\;\text{where $Q^{(2)}_{\mu,L}(x,v)$ as in \eqref{eq:Q_2}}. \label{eq:SO_finder}\\
				&\hspace{-0.5em} \text{Set }\;\; \alpha_k\eqdef\min \left\{1, 1/(2\norm{v^k}_{x^k})  \right\}. \label{eq:SO_alpha_k}
			\end{align}
					
			Set $z^{k}=x^{k} + \alpha_k v^k$, $i_k=i_k+1$\;
		}
		Set $M_{k+1} =\max\{\frac{L_{k}}{2},\underline{L}\}$, 
		$x^{k+1}=z^{k}$, $k=k+1$
	}
\end{algorithm}
\paragraph{Complexity bound.}
Our main result on the iteration complexity of Algorithm \ref{alg:SOAHBA} is the following Theorem. The proof follows similar steps as the proof of Theorem \ref{Th:AHBA_conv} and is given in Appendix \ref{sec:ProofSOM}.
\begin{theorem}
\label{Th:SOAHBA_conv}
Let Assumptions \ref{ass:1} and \ref{ass:2ndorder} hold. Fix the error tolerance $\eps>0$, the regularization parameter $\mu= \frac{\eps}{4\nu}$, and some initial guess $M_0>144\eps$ for the Lipschitz constant in \eqref{eq:SO_Lipschitz_Gradient}. Let $(x^{k})_{k\geq 0}$ be the trajectory generated by $\SAHBA(\mu,\eps,M_{0},x^{0})$, where $x^{0}$ is a $4\nu$-analytic center satisfying \eqref{eq:analytic_center}.
Then the algorithm stops in no more than 
\begin{equation}
\label{eq:SO_main_Th_compl}
\K_{II}(\eps,x^{0})= \ceil[\bigg]{\frac{576 \nu^{3/2} \sqrt{6\max\{M,M_0\}}(\Delta^f_0+ \eps)}{\eps^{3/2} }}
\end{equation}
outer iterations, and the number of inner iterations is no more than $2(\K_{II}(\eps,x^{0})+1)+2\max\{\log_{2}(2M/M_{0}),1\}$. Moreover, the output of $\SAHBA(\mu,\eps,M_{0},x^{0})$ constitutes an $(\eps,\frac{\max\{M,M_0\}\eps}{24\nu})$-2KKT point for problem \eqref{eq:Opt} in the sense of Definition \ref{def:eps_SOKKT}.
\end{theorem}
Note that Algorithm \ref{alg:SOAHBA} can be made anytime convergent by the same restarting procedure explained in Remark \ref{Rm:FO_restarts}.
\begin{remark}
\label{rem:SO_complexity_simplified}
Since $\Delta^f_0$ is expected to be larger than $\eps$, and the constant $M$ is potentially large, we see that the main term in the complexity bound \eqref{eq:SO_main_Th_compl} is $O\left(\frac{\nu^{3/2}\sqrt{M}\Delta^f_0}{\eps^{3/2}}\right)=O(\eps^{-3/2})$.
Note that the complexity result  $O(\max\{\eps_1^{-3/2},\eps_2^{-3/2}\})$ reported in \cite{CarDucHinSid19b,CarDucHinSid19} to find an $(\eps_1,\eps_2)$-2KKT point for arbitrary $ \eps_1,\eps_2 > 0$, is known to be optimal for unconstrained smooth non-convex optimization by second-order methods under the standard Lipschitz-Hessian assumption, subsumed on bounded sets by our Assumption \ref{ass:2ndorder}. A similar dependence on arbitrary $ \eps_1,\eps_2 > 0$ can be easily obtained from our theorem by setting $\eps=\min\{\eps_1,\eps_2\}$. 
\close
\end{remark}
\begin{remark}
An interesting observation is that our algorithm can be interpreted as a damped version of a cubic-regularized Newton's method. We have that the stepsize $\alpha_k$ satisfies $\alpha_k= \min \left\{1, \frac{1}{2\norm{v^k}_{x^k}}  \right\}$. At the initial phase, when the algorithm is far from an $(\eps_1,\eps_2)$-2KKT point, we have $\norm{v^k}_{x^k} > 1/2$ and $\alpha_k = \frac{1}{2\norm{v^k}_{x^k}}<1$. When the algorithm is getting closer to $(\eps_1,\eps_2)$-2KKT point, $\norm{v^k}_{x^k}$ becomes smaller and the algorithm automatically switches to full steps $\alpha_k=1$. 

At the same time our algorithm is completely different from cubic-regularized Newton's method \cite{NesPol06} applied to minimize the potential $F_{\mu}$. Indeed, we regularize by the cube of the local norm, rather than the cube of the standard Euclidean norm, and we do not form a second-order Taylor expansion of $F_{\mu}$. These adjustments are needed to align the search direction subproblem with the local geometry of the feasible set. Moreover, for our algorithm, the analysis of the cubic-regularized Newton's method is not directly applicable since it relies on stepsize 1, which may lead to infeasible iterates in our case.
\close
\end{remark}

\section{Conclusion}
In this paper, we propose two novel barrier algorithms for a large class of constrained non-convex optimization problems. Our algorithms possess favourable complexity guarantees that have the same dependence on target accuracy as the bounds for unconstrained optimization. Future works include the inexact implementation of our algorithms as well as applications to optimization problems in machine learning. 




\section*{Acknowledgements}
M. Staudigl acknowledges support from the COST Action CA16228 "European Network for Game Theory".
The work by P. Dvurechensky was supported by the Deutsche Forschungsgemeinschaft (DFG, German Research Foundation)
under Germany's Excellence Strategy – The Berlin Mathematics Research Center MATH$^+$ (EXC-2046/1, project ID: 390685689).




\bibliography{mybib}
\bibliographystyle{icml2024}

\newpage
\appendix
\onecolumn

\section{Discussion}
\label{sec:app_discussion}

In this section, we provide additional details to motivate our work and give concluding remarks. The following sections contain technical details of the proofs.

\subsection{Motivating Applications}


Consider a linear inverse problem with forward operator $\Phi$. Our aim is to learn a parameter $u\in U=\R^{n_{u}}_{+}$ so that the approximate equality 
$$
\Phi u\approx z 
$$
holds true. The matrix $\Phi$ has rows $\Phi_{1},\ldots,\Phi_{d}$ and we assume that $\Phi_{i}\in\R^{n_u}_{+}$. Furthermore, $z\in\R^{m}_{+}$ and $\sum_{j=1}^{p}\Phi_{ij}=\phi_{j}>0$ for all $j=1,\ldots,m$. In Poisson linear inverse problems \cite{HarMarWil11}, under additional structural assumptions, the penalized maximum likelihood approach for recovering the parameter $u$ leads to minimization of the function 
$$
\sum_{i=1}^{m}\{(\Phi u)_{i}-z_{i}\log((\Phi u)_{i})\}+\alpha r(u) 
$$
where $r(u)$ is a, potentially non-convex, regularizer. An important concrete example is the non-convex sparsity-inducing $\ell_{p}$ regularizer $r(u)=\sum_{i}\abs{u_{i}}^{p} = \norm{u}_p^p$, with $p\in(0,1)$, an example frequently used in computational statistics \cite{LohW2017}. We assume that $r$ is continuously differentiable on $\Int(U)$ and note that $\norm{u}_p^p$ is non-differentiable at $u=0$. 
This problem can be written in our optimization template by performing the variable substitution 
$
v=\Phi u 
$
and imposing the linear constraint 
$$
\Phi u-v=[\Phi;-I]\begin{pmatrix} u \\v \end{pmatrix}=0. 
$$
Therefore, defining the linear operator $\bA=[\Phi;-I]$, we obtain a matrix of rank $m$. We set $x=(u,v)$ and 
$$
f(x)=\sum_{i=1}^{m}\{v_{i}-z_{i}\log(v_{i})\}+\alpha r(u) 
$$
so that our inverse problem of Poisson image recovery reads as 
$$
\min_{x=(u,v)}f(x) \text{ s.t. } \bA x=0, x\in\bar{\setK}=\R^{n}_{+},
$$
where $n=n_{u}+m$. This problem admits the efficient self-concordant barrier 
$$
h(x)=-\sum_{i}\log(x_{i})\qquad\forall x\in\setK=\R^{n}_{++}. 
$$
By assumption, the function $f$ is twice continuously differentiable on $\R^{n}_{++}$. We compute the gradient and the Hessian as 
$$
\nabla f(x)=\begin{pmatrix} 
\alpha \nabla r(u) \\
\1_{m}-V^{-1}z
\end{pmatrix}, \text{ and } \nabla^{2}f(x)=\begin{pmatrix} \alpha \nabla^{2}r(u) & 0 \\
0 & V^{-2}z 
\end{pmatrix},
$$
where $V=\diag\{v_{1},\ldots,v_{n}\}$. The Hessian matrix of the barrier function decomposes as 
$$
H(x)=\nabla^{2}h(x)=\begin{pmatrix} U^{-2} & 0 \\0 & V^{-2}\end{pmatrix}.
$$
Thanks to the block structure of the involved matrices, the subproblems involved in the search direction finding routines of our two algorithms can be efficiently handled with efficient numerical linear algebra solvers. 

In a similar fashion, we can consider non-linear inverse problems where the loss is given by a squared misfit between the data $z$ and non-linear prediction function $\Phi(x)$ given, e.g., by a neural network. In this case, we have
\[
f(x)=\norm{\Phi(x)-z}^2+\alpha r(x). 
\]
The non-negativity constraints in this case may be motivated by the training of Input Convex Neural Networks \cite{pmlr-v70-amos17b}. A related problem of sparse non-linear regression may be reformulated as
\[
\min_{x \in \R^n} \norm{\Phi(x)-z}^2 \;\text {s.t.} \; \norm{x}_1\leq \lambda.
\]
This problem clearly fits our problem template with $f(x) = \norm{\Phi(x)-z}^2$ and $\setK = \{ x \in \R^n : \norm{x}_1\leq \lambda\}$.

\subsection{Additional Comments on Related Literature}
Interior-point methods \cite{DongDong2019Interior-Point} and optimization involving self-concordant functions \cite{Bac10,LinZha15,OstBac18,MarBacOstRu19,dvurechensky2020self-concordant,Carderera2021Simple} remains an active area of research in the Machine Learning community.  Yet, the main focus of this research stays on solving convex problems.  
However, the optimization community was recently quite successful in extending interior-point methods from the classical convex world \cite{NesIPM94} to non-convex world \cite{Ye92,TseBomSch11,HBA-linear,TseBomSch11,BiaCheYe15,HaeLiuYe18,NeiWr20,He:2022aa,dvurechensky2021hessian,Ye92,FayLu06,LuYua07}. The benefit of interior-point methods is that when the feasible set is given as an intersection of several sets, these methods allow decomposing the feasible set into separate building blocks. This allows one to avoid expensive projections onto the intersection. Further, such methods guarantee the feasibility of the iterates. 
At the same time, constrained optimization in the spirit of \cite{curtis2018complexity,hinder2018worst-case,cartis2019optimality,birgin2017complexity,grapiglia2020complexity,xie2019complexity} has recently attracted attention of ML community \cite{hong2023constrained} motivated by constrained deep neural networks, physical informed neural networks, PDE-constrained optimization, optimal control, and constrained model estimations.

With this paper, we further narrow the gap between the advances of optimization methods for non-convex problems with complicated constraints and Machine Learning applications. Moreover, our algorithms apply to more general problems than the ones available in the literature on interior-point methods for non-convex optimization, which we describe next and which influenced our work. In a sense, we generalize in this paper this line of works to a much more general class of problems.
The authors of \cite{HaeLiuYe18} propose first- and second-order algorithms with "optimal" \footnote{Here and below we refer to the complexity bound $O(\eps^{-2})$ for first-order and $O(\max\{\eps_1^{-3/2},\eps_2^{-3/2}\})$ for second-order methods as "optimal" for two reasons. First, the respective optimal algorithms for unconstrained problems have a similar dependence on the accuracy in the complexity bounds. Second, we are not aware of lower complexity bounds for problems with constraints that we consider. But, it is natural to expect that such problems are no easier than unconstrained problems.} complexity guarantees for problems with linear equality constraints and non-negativity constraints, i.e., $\bar{\setK}=\R_+^n$. Their algorithms are based on the Trust Region idea, unlike our algorithms that use quadratic and cubic regularization.
The authors of \cite{NeiWr20} consider a similar problem, but without equality constraints. They develop a Newton-Conjugate-Gradient method with "optimal" complexity to reach a second-order approximate stationary point. 
The authors of \cite{He:2022aa,dvurechensky2021hessian} consider a more general setting where $\bar{\setK}$ is a symmetric cone. \cite{dvurechensky2021hessian} propose first- and second-order methods with "optimal" complexity guarantees, and \cite{He:2022aa} propose a Newton-Conjugate-Gradient method with "optimal" complexity to reach a second-order approximate stationary point. 
Importantly, all these papers allow the objective to be non-differentiable at the relative boundary of the feasible set, unlike methods that use projections.

Compared to the optimization template \eqref{eq:Opt}, all these papers consider a narrower class of problems with $\bar{\setK}$ being a cone. Moreover, all the results in these papers heavily rely on the conic structure of the constraints (non-negativity constraints or general conic constraints). In particular, the conic structure allows one to easily introduce approximate optimality conditions since the conic duality can be used. Further, they use a narrow subclass of logarithmic or more generally logarithmically homogeneous self-concordant barriers that satisfy additionally
\[
h(tx)=h(x)-\nu\ln(t)\qquad \forall x\in\Int(\setK),t>0,
\]
and that are specific for cones and possess additional properties that can be used to derive optimality conditions and complexity results for algorithms. Our assumptions on the barrier are weaker, since we consider just self-concordant barriers without the additional logarithmic homogeneity property (which essentially forces the domain to be a pointed cone). This simultaneously allows us to solve much more general problems where $\bar{\setK}$ is a closed convex set, but not necessarily a cone. The latter in particular is justified by the existence of universal barriers for convex sets \cite{NesNem94}. We summarize the comparison with related literature in Table \ref{T:liter}.


\begin{table}[ht]
\label{T:liter}
\centering
\begin{tabular}{lcccccccc} 
\toprule
& Approach & Constraint & \makecell{ "Optimal"\\ complexity }  & No projection & \makecell{Feasible \\ iterates }   & Objective & Non-diff. \\ 
\midrule
\makecell{\cite{ghadimi2016mini-batch} \\ \cite{bogolubsky2016learning} \\ \cite{CarGouToi12} \\ \cite{birgin2017regularization} \\ \cite{CarGouToi18} } & Proximal &  Simple & $\surd$ & $\times$ & $\surd$ & General & $\times$\\
\midrule
\makecell{\cite{Andreani:2019uf} \\ \cite{Andreani:2021vq}  } & \makecell{ Augmented \\ Lagrangian  } &  Only cones & $\times$ & $\surd$ & $\times$ & General & $\times$\\
\midrule
\cite{LuYua07} & Interior-point &  Only cones & $\times$ & $\surd$ & $\surd$ & Quadratic & $\times$\\
\midrule
\cite{BiaCheYe15} & Interior-point &  Box & $\surd$ & $\surd$ & $\surd$ & Structured & $\surd$\\
\midrule
\makecell{\cite{TseBomSch11} \\ \cite{HBA-linear}  } &Interior-point &  Only $\R_+^n$ & ? & $\surd$ & $\surd$ & Quadratic & $\times$\\
\midrule
\makecell{ \cite{HaeLiuYe18} \\ \cite{NeiWr20} } &Interior-point &  Only $\R_+^n$ & $\surd$ & $\surd$ & $\surd$ & General & $\surd$\\
\midrule
\makecell{\cite{dvurechensky2021hessian} \\ \cite{He:2022aa} } & Interior-point &  Only cones & $\surd$ & $\surd$ & $\surd$ & General & $\surd$ \\
\midrule
This paper   & Interior-point & \textbf{General} 
& $\surd$ & $\surd$ & $\surd$ & General & $\surd$ \\ 
\bottomrule
\end{tabular}
\caption{Summary of literature. "No projection" means that the algorithm does not need to project onto the complicated $\bar{\setX}$. "Non-diff." means that the objective may be non-differentiable at the relative boundary of the feasible set.}
\end{table}

\subsection{Conclusion}

In this paper, we propose first- and second-order alrgorithms for non-convex problems with linear and general set constraints. We  develop also necessary optimality conditions for such problems and define their suitable approximate counterparts. Further, we show that our algorithms achieve approximate stationary points with "optimal" worst-case iteration complexity. Unlike previously known results on interior-point methods for non-convex optimization, our approach allows one to solve a much wider class of problems.

Future works include extensions of our algorithms for the setting of inexact solution of search direction finding problems. With that respect we believe that it is possible to construct a Newton-Conjugate-Gradients counterpart of out second-order method. Further, we plan to use the proposed methods in machine learning applications such as constrained non-linear regression and training input convex neural networks. Further potential extensions include adding non-linear functional constraints to the problem.

\section{Auxiliary Facts on SCBs}
\label{app:barrier}

The following properties are taken from \cite{Nes18}, Lemma 5.4.3, Theorem 5.3.7.
\begin{proposition}
\label{prop:logSCB}
Let $h\in\scrH_{\nu}(\setK)$, $x\in\setK$, $t>0$ and $H(x)=\nabla^{2}h(x)$.
Then,
\begin{align}
&\inner{ \nabla h(x), [H(x)]^{-1}\nabla h(x) } \leq \nu. \label{eq:log_hom_scb_norm_prop}\\
&\inner{\nabla h(x),y-x}<\nu\qquad\forall y\in\setK.\label{eq:SCBangle}
\end{align}
\end{proposition}
Note that \eqref{eq:SCBangle} means that $\nabla h(x) \in \NC_{\bar{\setK}}^{\nu}(x)$. 
The following fact may be derived from \eqref{eq:SCBangle} and is its appriximate counterpart.
\begin{proposition}[{Proposition 2.7 \cite{MonSicSva15}}]
\label{prop:SCB_normal_cone}
Let $h\in\scrH_{\nu}(\setK)$, $x\in\setK$, and $s \in \setE$ satisfy $\norm{s-\nabla h(x)}_{x}^*\leq \xi <1$. Then, for all $y \in \bar{\setK}$,
\begin{equation}
\label{eq:approx_normal_cone}
\inner{s,y-x} \leq \nu + \frac{\sqrt{\nu}+\xi}{1-\xi} \xi.
\end{equation}
\end{proposition}

\section{Proof of Theorem \ref{Th:limit_optimality}}
\label{sec:optimality_proof}

Let  $x^*$ be a local solution for problem \eqref{eq:Opt}.
We consider the following perturbed version of problem \eqref{eq:Opt}, for which $x^*$ is the unique global solution when $\delta>0$ sufficiently small, 
\begin{equation}
\label{eq:Haeser_limit_optimality_proof_1}
	\min_{x} f(x) + \frac{1}{4}\|x-x^*\|^4 \quad \text{s.t.: } \bA x= b,\; x \in \bar{\setK}, \; \|x-x^*\|^2 \leq \delta.
\end{equation}
Next, using the barrier $h$ for $\setK$, we change the constraint $x\in\setK$ to the penalty $\mu_k h(x)$, where $\mu_k> 0$, $\mu_k \to 0$ is a given sequence. This leads us to the following parametric sequence of problems for $k\geq 0$
\begin{equation}
\label{eq:Haeser_limit_optimality_proof_2}
\begin{array}{ll}
	&\min_{x} f_k(x) \eqdef f(x) + \frac{1}{4}\|x-x^*\|^4 + \mu_k h(x)\\
\text{s.t.: }& \bA x= b,\; x\in\Int(\setK), \; \|x-x^*\|^2 \leq \delta.
\end{array}
\end{equation}
From the classical theory of interior penalty methods \cite{FiacMcCo68}, it is known that a global solution $x^k$ exists for this problem for all $k$ and that cluster points of $x^k$ are global solutions of \eqref{eq:Haeser_limit_optimality_proof_1}. Clearly, $ \bA x^k = b$ and $x^k\in\setK=\Int(\setK)=\Int(\bar{\setK})$. Since $\|x^k-x^*\|^2 \leq \delta$, the sequence $x^k$ is bounded and, thus, $x^k \to x^*$. This finishes the proof of item \ref{Th:limit_optimality_1} of Theorem \ref{Th:limit_optimality}.

Since $x^k \to x^*$, for large enough $k$, $x^k$ is a local solution to the problem 
\begin{equation}
\label{eq:Haeser_limit_optimality_proof_3}
	\min_{x}  f(x) + \frac{1}{4}\|x-x^*\|^4 + \mu_k h(x)  \quad  \text{s.t.: } \; \bA x = b
\end{equation}
and $x^k \in \setK$.
By Assumption \ref{ass:1}, the system of constraints $ \bA x = b$ has full rank. Hence, we can write necessary optimality conditions for problem \eqref{eq:Haeser_limit_optimality_proof_3} which say that there exists a Lagrange multiplier $y^k \in \R^m$ such that
\begin{equation}
\label{eq:Haeser_limit_optimality_proof_4}
0 = \nabla f(x^k) + \|x^k-x^*\|^2(x^k-x^*) - \bA y^k  +   \mu_k \nabla h(x^k).
\end{equation}
Let us choose the vector $s^k = -\mu_k \nabla h(x^k)$.  Then, item \ref{Th:limit_optimality_2} of Theorem \ref{Th:limit_optimality} follows from \eqref{eq:Haeser_limit_optimality_proof_4} since $\|x^k-x^*\|^2(x^k-x^*) \to 0$. 
Further, by \eqref{eq:SCBangle} with $x=x^{k}$, we have for any $y \in \bar{\setK}$
\begin{align*}
\inner{\nabla h(x^{k}),y-x^{k}} \leq \nu \Leftrightarrow \inner{-s^k/\mu_k,y-x^{k}} \leq \nu \Leftrightarrow \inner{-s^k,y-x^{k}} \leq \mu_k\nu. 
\end{align*}
Taking $\sigma_k=\mu_k\nu \to 0$, we see that $-s^k \in \NC_{\bar{\setK}}^{\sigma_k}(x^{k})$. This proves item \ref{Th:limit_optimality_3} of Theorem \ref{Th:limit_optimality}. 

The second-order differentiability assumption and the full rank condition give the following second-order necessary optimality condition for \eqref{eq:Haeser_limit_optimality_proof_3}. For all $d \in \setE$ such that $\bA d=0$, it holds that
\begin{equation}
\inner{(\nabla^2 f(x^k) + \mu_k  H(x^k) +\Sigma_{k})d,d}\geq 0,
 \end{equation}
where $\Sigma_k = 2(x^k-x^*)(x^k-x^*)^{\top} +   \|x^k-x^*\|^{2}\bI$, $\bI$ being the identity operator. Setting $\theta_k=\mu_k$  and $\delta_k $ as the largest eigenvalue of the positive semi-definite matrix $\Sigma_k$, we conclude that $\theta_{k}\to 0$ and $\delta_{k}\to 0$ as $k\to\infty$. This finishes the proof of eq. \eqref{Th:limit_optimality_4} and Theorem \ref{Th:limit_optimality}. \qed

\section{Missing Proofs from Section \ref{sec:firstorder}}
\label{sec:App_FO}

\subsection{Proof of Theorem \ref{Th:AHBA_conv}}
\label{sec:ProofFOM}
Our proof proceeds in four steps. First, we show that $\AHBA(\mu,\eps,L_{0},x^{0})$ produces points in $\setX$, and, thus, is indeed an interior-point method. Then, we proceed to show that the line-search process of finding appropriate $L_k$'s in each iteration is finite, and estimate the total number of attempts in this process. After that, we prove that if the stopping criterion does not hold at iteration $k$, i.e., $\norm{v^{k}}_{x^{k}} \geq \tfrac{\eps}{3\nu}$, then the objective $f$ is decreased by a quantity $O(\eps^{2})$. From the global lower boundedness of the objective, we derive that the method stops in at most $O(\eps^{-2})$ iterations. Finally, we show the opposite, i.e., that when the stopping criterion holds, the method has generated an $\eps$-KKT point according to Definition \ref{def:eps_KKT}. 

\subsubsection{Interior-Point Property of the Iterates}
\label{S:FO_correct}
We start the induction argument by observing that $x^{0}\in\setX$ by construction. Further, let $x^{k}\in\setX = \setK \cap \setL$ be the $k$-th iterate of the algorithm with the corresponding step direction $v^{k}\eqdef v_{\mu}(x^{k})$. 
By eq. \eqref{eq:alpha_k}, the step-size $\alpha_k$ satisfies $\alpha_{k}\leq \frac{1}{2\norm{v^k}_{x^k}}$, and, hence, $\alpha_{k}\norm{v^k}_{x^k}\leq 1/2$ for all $k\geq 0$. 
Thus, by Lemma \ref{lem:Dikin}, we have $x^{k+1}=x^{k}+\alpha_{k}v^{k} \in\setK$. Since, by \eqref{eq:finder}, $\bA v^{k} =0$, we have that $x^{k+1} \in \setL$. Thus, $x^{k+1}\in\setK \cap \setL=\setX$. By induction, we conclude that $(x^{k})_{k\geq 0}\subset\setX$. 

\subsubsection{Bounding the Number of Backtracking Steps}
\label{sec:backtrack1}
Consider iteration $k$. The sequence $2^{i_k} L_k $ is increasing as $i_k$ is increasing. Hence, by Assumption \ref{ass:gradLip}, we know that when $2^{i_k} L_k \geq \max\{M,L_k\}$, the line-search process for sure stops since inequality \eqref{eq:LS} holds. 
Thus, $2^{i_k} L_k \leq 2\max\{M,L_k\}$ must be the case, and, consequently, $L_{k+1} = 2^{i_k-1} L_k \leq \max\{M,L_k\}$, which, by induction, gives $L_{k+1} \leq \bar{M}\eqdef\max\{M,L_0\}$. 
At the same time, $\log_{2}\left(\frac{L_{k+1}}{L_{k}}\right)= i_{k}-1$, $\forall k\geq 0$. Let $N(k)$ denote the number of inner line-search iterations up to the $k-$th iteration of $\AHBA(\mu,\eps,L_{0},x^{0})$. Then, using that $L_{k+1} \leq \bar{M}=\max\{M,L_0\}$, we obtain
\begin{align*}
N(k)&=\sum_{j=0}^{k}(i_{j}+1)=\sum_{j=0}^{k} (\log_{2}(L_{j+1}/L_{j})+2 ) \leq 2(k+1)+\max\{\log_{2}(M/L_{0}),0\}.
\end{align*}
As we see, on average the inner loop ends after two trials. 

\subsubsection{Bound for the Number of Outer Iterations}
Our goal now is to establish per-iteration decrease of the potential $F_{\mu}$ if the stopping condition does not yet hold, i.e., $\norm{v^{k}}_{x^{k}} \geq \tfrac{\eps}{3\nu}$, and derive from that a global iteration complexity bound for $\AHBA$.
Let us fix iteration counter $k$. Since $L_{k+1} = 2^{i_k-1}L_k$, the step-size \eqref{eq:alpha_k} is equivalent to $\alpha_{k}=\min \left\{\frac{1}{2L_{k+1} + 2 \mu},\frac{1}{2\norm{v^k}_{x^{k}}} \right\}$. Hence, $\alpha_{k}\norm{v^k}_{x^{k}}\leq 1/2$, and \eqref{eq:success} with the  substitution $t=\alpha_k = \ct_{\mu,2L_{k+1}}(x^{k})$, $M=2L_{k+1}$, $x=x^{k}$, $v_{\mu}(x^{k}) \eqdef v^k$ gives:  
\begin{equation}
\label{eq:FO_per_iter_proof_2}
F_{\mu}(x^{k+1})-F_{\mu}(x^{k})\leq -\alpha_k \norm{v^k}_{x^k}^{2}\left(1-(L_{k+1}+\mu)\alpha_k\right) \leq -\frac{\alpha_k \norm{v^k}_{x^k}^{2}}{2},
\end{equation}
where the last inequality is since $\alpha_k \leq \frac{1}{2(L_{k+1}+\mu)}$. 
Substituting into \eqref{eq:FO_per_iter_proof_2} the two possible values of the step-size $\alpha_k$ in \eqref{eq:alpha_k} gives
\begin{equation}
\label{eq:per_iter_decr_0}
F_{\mu}(x^{k+1})-F_{\mu}(x^{k})\leq 
\left\{
\begin{array}{ll}
- \frac{\norm{v^{k}}_{x^{k}}^{2} }{4(L_{k+1}+\mu)} & \text{if}  \alpha_k=\frac{1}{2(L_{k+1}+\mu)},\\

 - \frac{\norm{v^{k}}_{x^{k}}}{4} & \text{if }  \alpha_k=\frac{1}{2\norm{v^k}_{x^{k}}}.
\end{array}\right.
\end{equation}
As we proved in Section \ref{sec:backtrack1}, $L_{k+1} \leq \bar{M}$. Thus, we obtain that 
\begin{equation}
\label{eq:per_iter_decr}
F_{\mu}(x^{k+1}) - F_{\mu}(x^{k}) \leq  -\frac{\norm{v^{k}}_{x^{k}}}{4} \min\left\{1,  \frac{\norm{v^{k}}_{x^{k}}}{\bar{M}+\mu}\right\} \eqdef -\delta_{k}.
\end{equation}
Rearranging and summing these inequalities for $k$ from $0$ to $K-1$ gives
\begin{align}
&K\min_{k=0\ldots,K-1} \delta_{k} \leq  \sum_{k=0}^{K-1}\delta_{k}\leq F_{\mu}(x^{0})-F_{\mu}(x^{K}) \notag \\
& \quad\stackrel{\eqref{eq:potential}}{=} f(x^{0}) - f(x^{K}) + \mu (h(x^{0}) - h(x^{K})) \leq f(x^{0}) - f_{\min}(\setX) + \eps, \label{eq:FO_per_iter_proof_6} 
\end{align}
where in the last inequality we used that, by the assumptions of Theorem \ref{Th:AHBA_conv}, $x^{0}$ is a $\nu$-analytic center defined in \eqref{eq:analytic_center} and $\mu = \eps/\nu$, implying that $h(x^{0}) - h(x^{K}) \leq \nu = \eps/\mu$.
Thus, up to passing to a subsequence, $\delta_{k}\to 0$, and consequently $\norm{v^{k}}_{x^{k}} \to 0$ as $k \to \infty$. Hence, the stopping criterion in Algorithm \ref{alg:AHBA} is achievable and the algorithm is correctly defined in this respect.

Let us now assume that the stopping criterion $\norm{v^k}_{x^k} < \frac{\eps}{3\nu}$ does not hold for $K$ iterations of $\AHBA$. Then, for all $k=0,\ldots,K-1,$ we have 
$\delta_{k}\geq \min\left\{\frac{\eps}{12 \nu},\frac{\eps^{2}}{36 \nu^{2}(\bar{M}+\mu)}\right\}$. 
Using that we set $\mu=\frac{\eps}{\nu}$, it follows from \eqref{eq:FO_per_iter_proof_6} that
\[
K\frac{\eps^{2}}{36 \nu^{2}(\bar{M}+\eps/\nu)} = K \min\left\{\frac{\eps}{12 \nu},\frac{\eps^{2}}{36 \nu^{2}(\bar{M}+\eps/\nu)}\right\}\leq f(x^{0})-f_{\min}(\setX)+\eps.
\]
Hence, recalling that $\bar{M}=\max\{M,L_0\}$, we obtain
\[
K \leq 36(f(x^{0}) - f_{\min}(\setX)+ \eps) \cdot \frac{\nu^{2}(\max\{M,L_0\}+\eps/\nu)}{\eps^{2}}.
\] 
Thus, we obtain the bound on the number of iterations on which the stopping criterion is not satisfied.
This, combined with the bound for the number of inner steps in Section \ref{sec:backtrack1}, proves the complexity bound statement of Theorem \ref{Th:AHBA_conv}.

\subsubsection{Generating $\eps$-KKT Point}
To finish the proof of Theorem \ref{Th:AHBA_conv}, we now show that when Algorithm \ref{alg:AHBA} stops for the first time, it returns a $2\eps$-KKT point of \eqref{eq:Opt} according to Definition \ref{def:eps_KKT}.

Clearly, \eqref{eq:eps_optim_equality_set} in Definition \ref{def:eps_KKT} holds by the construction of the algorithm. Thus, we focus on showing \eqref{eq:eps_optim_grad}.
Let the stopping criterion hold at iteration $k$, i.e., $\norm{v^{k}}_{x^{k}} < \frac{\eps}{3\nu}$. Using the optimality condition \eqref{eq:finder_1} at iteration $k$ and the definition of the potential \eqref{eq:potential}, we get
\begin{equation}
\label{eq:FO_KKT_proof_0}
\nabla f(x^{k})-\bA^{\ast}y^{k}+\mu \nabla h(x^{k}) =-H(x^{k})v^{k}\iff [H(x^{k})]^{-1}\left(\nabla f(x^{k})-\bA^{\ast}y^{k}+\mu \nabla h(x^{k}) \right)=-v^{k}.
\end{equation}

Multiplying both equations, taking the square root, and using the stopping criterion $\norm{v^{k}}_{x^{k}} < \frac{\eps}{3\nu}$ 
we obtain 
\begin{equation}
\label{eq:FO_KKT_proof_00}
\norm{\nabla f(x^{k})-\bA^{\ast}y^{k}+\mu \nabla h(x^{k})}^{\ast}_{x^{k}}=\norm{v^{k}}_{x^{k}}<\frac{\eps}{3\nu},
\end{equation}
whence, dividing by $\mu$
\begin{equation}
\label{eq:FO_KKT_proof_01}
\norm{-\frac{1}{\mu}(\nabla f(x^{k})-\bA^{\ast}y^{k})-\nabla h(x^{k})}^{\ast}_{x^{k}}<\frac{\eps}{3\mu\nu}=\frac{1}{3}.
\end{equation}
Applying \eqref{eq:approx_normal_cone} with $s=-\frac{1}{\mu}(\nabla f(x^{k})-\bA^{\ast}y^{k})$ and  $\xi=\frac{\eps}{3\mu\nu}=\frac{1}{3}$, we obtain
\begin{equation}
\label{eq:FO_KKT_proof_02}
\inner{-\frac{1}{\mu}(\nabla f(x^{k})-\bA^{\ast}y^{k}),x-x_k}<\nu+\frac{\sqrt{\nu}+\xi}{1-\xi} \xi = \nu + \frac{\sqrt{\nu}+1/3}{2}\qquad \forall x \in \bar{\setK},
\end{equation}
whence,  
\begin{equation}
\label{eq:FO_KKT_proof_03}
\inner{\nabla f(x^{k})-\bA^{\ast}y^{k},x-x_k}>-\mu\nu-\mu\frac{\sqrt{\nu}+1/3}{2} > -2\eps \qquad \forall x \in \bar{\setK},
\end{equation}
where we used that $\mu=\frac{\eps}{\nu}$ and that $\nu\geq 1$. Thus, we obtain that \eqref{eq:eps_optim_grad} holds, which finishes the proof of Theorem \ref{Th:AHBA_conv}.

\section{Missing Proofs from Section \ref{sec:secondorder}}
\label{sec:App_SO}

\subsection{Proofs of Preliminary Results}

\paragraph{Proof of the implication \eqref{eq:LipHess} $\Rightarrow$ \eqref{eq:SO_Lipschitz_Gradient}.}
We have
\begin{align*}
&\norm{\nabla f(x+v)-\nabla f(x)-\nabla^{2}f(x)v}^{\ast}_{x}=\norm{\int_{0}^{1}(\nabla^{2}f(x+tv)-\nabla^{2}f(x))v\dif t }^{\ast}_{x}\\
&\leq \int_{0}^{1}\norm{\nabla^{2}f(x+tv)-\nabla^{2}f(x)}_{\text{op},x}\cdot \norm{v}_{x}\dif t   \leq \frac{M}{2}\norm{v}^{2}_{x}.
\end{align*}

\paragraph{Proof of \eqref{eq:cubicestimate}.}
To obtain \eqref{eq:cubicestimate}, observe that for all $x\in\setX$ and $v\in\scrT_{x}$, we have
\begin{align*}
&\abs{f(x+v)-f(x)-\inner{\nabla f(x),v}-\frac{1}{2}\inner{\nabla^{2}f(x)v,v}}=\abs{\int_{0}^{1}\inner{\nabla f(x+tv)-\nabla f(x)-\frac{1}{2}\nabla^{2}f(x)v,v}\dif t}\\
&\quad\leq \int_{0}^{1}\norm{\nabla f(x+tv)-\nabla f(x)-\frac{1}{2}\nabla^{2}f(x)v}_{x}^{*}\dif t\cdot\norm{v}_{x} 
\leq \frac{M}{6}\norm{v}^{3}_{x}.
\end{align*}

\paragraph{Proof of Proposition \ref{Th:PD}.}
The high-level idea is that after a transition to the basis induced by the affine subspace $\setL_{0}$, the subproblem \eqref{eq:cubicproblem} becomes an unconstrained minimization problem similar to the Cubic Newton step in \cite{NesPol06}.
To that end, let $\{z_{1},\ldots,z_{p}\}$ be an orthonormal basis of $\setL_{0}$ and the linear operator $\bZ:\R^{p}\to\setL_{0}$ be defined by $\bZ w=\sum_{i=1}^{p}z_{i}w^{i}$ for all $w=[w^{1};\ldots;w^{p}]^{\top}\in\R^{p}$.
Based on this linear map, we define the projected data
\begin{equation}\label{eq:dataKernel}
\gbold\eqdef\bZ^{\ast}\nabla F_{\mu}(x),\; \bJ\eqdef\bZ^{\ast}\nabla^{2}f(x)\bZ,\;\bH\eqdef\bZ^{\ast}H(x)\bZ \succ 0
\end{equation}
and apply it, together with the change of variables $v=\bZ u$, to reformulate the search-direction finding problem \eqref{eq:cubicproblem} as an unconstrained cubic-regularized subproblem of finding $u_{L}\in\R^{p}$ s.t.
\begin{equation}\label{eq:cubicauxiliary}
u_{L}\in  \Argmin_{u\in\R^{p}}\{\inner{\gbold,u}+\frac{1}{2}\inner{\bJ u,u}+\frac{L}{6}\norm{u}^{3}_{\bH}\},
\end{equation}
where  $\norm{\cdot}_{\bH}$ is the norm induced by the operator $\bH$. From \cite{NesPol06}, Thm. 10 we deduce  
\[
\bJ+\frac{L\norm{u_{L}}_{\bH}}{2}\bH\succeq 0.
\]
Denoting $v_{\mu,L}(x) = \bZ u_{L}$, we see
\begin{align*}
\norm{u_{L}}_{\bH}=\inner{\bZ^{\ast}H(x)\bZ u_{L},u_{L}}^{1/2}=\inner{H(x)(\bZ u_{L}),\bZ u_{L}}^{1/2}&=\norm{v_{\mu,L}(x)}_{x}, \text{  and} \\
\bZ^{\ast}\left(\nabla^{2}f(x)+\frac{L}{2}\norm{v_{\mu,L}(x)}_{x}H(x)\right)\bZ&\succeq 0,
\end{align*}
which implies $\nabla^{2}f(x)+\frac{L}{2}\norm{v_{\mu,L}(x)}_{x}H(x)\succeq 0$ over the nullspace $\setL_{0} = \{ v\in\setE:\bA v=0\}$. 
\qed

\noindent
The above derivations give us a hint on possible approaches to numerically solve problem \eqref{eq:cubicproblem} in practice.
Before the start of the algorithm, as a preprocessing step, we once calculate matrix $\bZ$ and use it during the whole algorithm execution. At each iteration, we calculate the new data using \eqref{eq:dataKernel} and get a standard \textit{unconstrained} cubic subproblem  \eqref{eq:cubicauxiliary}. \cite{NesPol06} show how such problems can be transformed to a \emph{convex} problem to which fast convex programming methods could in principle be applied. However, we can also solve it via recent efficient methods based on Lanczos' method \cite{CarGouToi11a,Jia22}. In any case, we can recover our step direction $v_{\mu,L}(x)$ by the matrix vector product $\bZ u_{L}$, where  $u_{L}$ is the solution obtained from this subroutine.

\paragraph{Derivation of the step-size of Algorithm \ref{alg:SOAHBA} and derivation of \eqref{eq:progressSOM_main}.} Our goal now is to construct an admissible step-size policy,  given the step direction $v_{\mu,L}(x)$. We act in a similar fashion as in the analysis of the first-order algorithm by applying optimality conditions and estimating the per-iteration decrease of the potential depending on the step-size.
Let $x\in\setX$ be the current position of the algorithm. Define  $x^{+}(t)\eqdef x+t v_{\mu,L}(x)$, where $t\geq 0$ is a step-size. 
By Lemma \ref{lem:Dikin} and since $v_{\mu,L}(x)\in\setL_{0}$ by \eqref{eq:opt2}, we know that $x^{+}(t)$ is in $\setX$ provided that $t\in I_{x,\mu,L}\eqdef[0,\frac{1}{\norm{v_{\mu,L}(x)}_x})$. For all such $t$, by \eqref{eq:cubicdecrease}, we get 
\begin{equation}\label{eq:SO_potent_upper_bound}
\begin{split}
F_{\mu}(x^{+}(t))\leq F_{\mu}(x)&+t\inner{\nabla F_{\mu}(x),v_{\mu,L}(x)} + \frac{t^2}{2}\inner{\nabla^{2}f(x)v_{\mu,L}(x),v_{\mu,L}(x)} \\
&+ \frac{Mt^3}{6}\norm{v_{\mu,L}(x)}^{3}_{x}  +\mu t^{2}\norm{v_{\mu,L}(x)}_x^2\omega(t\norm{v_{\mu,L}(x)}_x).
\end{split}
\end{equation}
Since $v_{\mu,L}(x) \in \setL_0 = \{ v \in \setE \vert \bA v = 0\}$, multiplying \eqref{eq:PD} with $v_{\mu,L}(x)$ from the left and the right, and multiplying \eqref{eq:opt1} by $v_{\mu,L}(x)$ and combining with \eqref{eq:opt2}, we obtain
\begin{align}
&\inner{\nabla^{2}f(x)v_{\mu,L}(x),v_{\mu,L}(x)}\geq-\frac{L}{2}\norm{v_{\mu,L}(x)}^{3}_{x},\label{eq:descent1}\\
&\inner{\nabla F_{\mu}(x),v_{\mu,L}(x)}+\inner{\nabla^{2}f(x)v_{\mu,L}(x),v_{\mu,L}(x)}+\frac{L}{2}\norm{v_{\mu,L}(x)}^{3}_{x}=0.\label{eq:normal}
\end{align}
Under the additional assumption that $t\leq 2$ and $L\geq M$, we obtain
\begin{align*}
&t\inner{\nabla F_{\mu}(x),v_{\mu,L}(x)} + \frac{t^2}{2}\inner{\nabla^{2}f(x)v_{\mu,L}(x),v_{\mu,L}(x)} + \frac{Mt^3}{6}\norm{v_{\mu,L}(x)}^{3}_{x}\\
&\stackrel{\eqref{eq:normal}}{=} - t \left(\inner{\nabla^{2}f(x)v_{\mu,L}(x),v_{\mu,L}(x)}+\frac{L}{2}\norm{v_{\mu,L}(x)}^{3}_{x}\right) &\\
&+ \frac{t^2}{2}\inner{\nabla^{2}f(x)v_{\mu,L}(x),v_{\mu,L}(x)} + \frac{Mt^3}{6}\norm{v_{\mu,L}(x)}^{3}_{x} \\
& = \left(\frac{t^2}{2} - t \right) \inner{\nabla^{2}f(x)v_{\mu,L}(x),v_{\mu,L}(x)} - \frac{Lt}{2}\norm{v_{\mu,L}(x)}^{3}_{x} + \frac{Mt^3}{6}\norm{v_{\mu,L}(x)}^{3}_{x} \\
& \stackrel{\eqref{eq:descent1},t\leq 2}{\leq} \left(\frac{t^2}{2} - t \right) \left(- \frac{L}{2}\norm{v_{\mu,L}(x)}^{3}_{x} \right) - \frac{Lt}{2}\norm{v_{\mu,L}(x)}^{3}_{x} + \frac{Mt^3}{6}\norm{v_{\mu,L}(x)}^{3}_{x} \\
& = - \norm{v_{\mu,L}(x)}^{3}_{x} \left(\frac{Lt^2}{4} - \frac{Mt^3}{6} \right) 
\stackrel{L \geq M}{\leq} - \norm{v_{\mu,L}(x)}^{3}_{x} \frac{Lt^2}{12} \left(3 - 2t \right).
\end{align*}
Substituting this into \eqref{eq:SO_potent_upper_bound}, we arrive at
\begin{align*}
F_{\mu}(x^{+}(t))&\leq F_{\mu}(x) - \norm{v_{\mu,L}(x)}^{3}_{x} \frac{Lt^2}{12} \left(3 - 2t \right)+\mu t^{2}\norm{v_{\mu,L}(x)}_x^2\omega(t\norm{v_{\mu,L}(x)}_x)\\ 
&\stackrel{\eqref{eq:omega_upper_bound}}{\leq}  F_{\mu}(x) - \norm{v_{\mu,L}(x)}^{3}_{x} \frac{Lt^2}{12} \left(3 - 2t \right)+\mu \frac{t^{2}\norm{v_{\mu,L}(x)}_{x}^{2}}{2(1-t\norm{v_{\mu,L}(x)}_x)}. 
\end{align*} 
for all  $t \in I_{x,\mu,L}$. Therefore, if $t\norm{v_{\mu,L}(x)}_x\leq 1/2$, we finally obtain
\begin{align}
F_{\mu}(x^{+}(t))-F_{\mu}(x)&\leq - \frac{Lt^2\norm{v_{\mu,L}(x)}^{3}_{x} }{12} \left(3 - 2t \right)+\mu t^{2}\norm{v_{\mu,L}(x)}_{x}^{2} \nonumber \\
&= - \norm{v_{\mu,L}(x)}^{3}_{x} \frac{Lt^2}{12}\left(3 - 2t - \frac{12 \mu}{L\norm{v_{\mu,L}(x)}_{x}} \right) \eqdef -\eta_{x}(t).\label{eq:progressSOM}
\end{align}
This is exactly the bound \eqref{eq:progressSOM_main}.
Unfortunately, finding and using the explicit maximizer of $\eta_{x}(t)$ is quite challenging. But, as we will see, the following step-size is a good and simple alternative:
\begin{equation}
\label{eq:SO_t_opt_def}
\ct_{\mu,L}(x) \eqdef \frac{1}{\max\{1,2\norm{v_{\mu,L}(x)}_x\}} = \min \left\{1, \frac{1}{2\norm{v_{\mu,L}(x)}_x}  \right\}.
\end{equation}  
Note that $\ct_{\mu,L}(x) \leq 1$ and $\ct_{\mu,L}(x)\norm{v_{\mu,L}(x)}_x\leq 1/2$. Thus, this choice of the step-size is feasible to derive \eqref{eq:progressSOM}.

\subsection{Proof of Theorem \ref{Th:SOAHBA_conv}}
\label{sec:ProofSOM}
The main steps of the proof are similar to the analysis of Algorithm \ref{alg:AHBA}. We start by showing the feasibility of the iterates and the correctness of the backtracking line-search process, i.e., that this process is finite. We also estimate the total number of attempts in this process. After that, we analyze the per-iteration decrease of $F_{\mu}$ and show that if the stopping criterion does not hold at iteration $k$, then the objective function is decreased by the value $O(\eps^{3/2})$. This, by the global lower boundedness of the objective, allows us to conclude that the algorithm stops in  $O(\eps^{-3/2})$ iterations.
 Finally, we show the opposite, i.e., that when the stopping criterion holds, the method has generated an approximate second-order KKT point in the sense of Definition \ref{def:eps_SOKKT}.

\subsubsection{Interior-Point Property of the Iterates}
We start the induction argument by observing that $x^{0}\in\setX$ by construction. Further, let $x^{k}\in\setX = \setK \cap \setL$ be the $k$-th iterate of the algorithm with the corresponding step direction $v^{k}\eqdef v_{\mu,L}(x^{k})$.  By  \eqref{eq:SO_alpha_k}, the step-size $\alpha_k$ satisfies $\alpha_{k}\leq \frac{1}{2\norm{v^{k}}_{x^{k}}}$. Consequently, $\alpha_{k}\norm{v^{k}}_{x^{k}}\leq 1/2$ for all $k\geq 0$, and using Lemma \ref{lem:Dikin} as well as equality $\bA v^{k} =0$ by \eqref{eq:SO_finder}, we have that $x^{k+1}=x^{k}+\alpha_{k}v^{k}\in\setK \cap \setL = \setX$. By induction, it follows that $x^{k}\in\setX$ for all $k\geq 0$.

\subsubsection{Bounding the Number of Backtracking Steps}
\label{sec:backtrack2}
To bound the number of cycles involved in the line-search process for finding appropriate constants $L_{k}$, we proceed as in Section \ref{sec:backtrack1}. 
Let us fix an iteration $k$. The sequence $L_k = 2^{i_k} M_k$ is increasing as $i_k$ is increasing, and Assumption \ref{ass:2ndorder} holds. 
This implies \eqref{eq:cubicestimate}, and thus when $L_k = 2^{i_k} M_k \geq \max\{M,M_k\}$, the line-search process for sure stops since inequalities \eqref{eq:SO_LS_1} and \eqref{eq:SO_LS_2} hold.	
Hence, $L_k=2^{i_k} M_k \leq 2\max\{M,M_k\}$ must be the case, and, consequently, $M_{k+1}=\max\{L_k/2, \underline{L}\} \leq  \max\{\max\{M,M_k\}, \underline{L}\} = \max\{M,M_k\}$, which, by induction, gives $M_{k} \leq \bar{M} \eqdef \max\{M,M_0\}$ and $L_k  \leq 2\bar{M}$. 
%
%
%
%
At the same time, by construction, $M_{k+1}= \max\{2^{i_k-1}M_{k},\underline{L}\} = \max\{L_k/2,\underline{L}\} \geq L_k/2 $. Hence, $L_{k+1} = 2^{i_{k+1}} M_{k+1} \geq 2^{i_{k+1}-1} L_k$ and therefore $\log_{2}\left(\frac{L_{k+1}}{L_{k}}\right)\geq i_{k+1}-1$, $\forall k\geq 0$. At the same time, at iteration $0$ we have $L_0=2^{i_0} M_0 \leq 2\bar{M}$, whence, $i_0 \leq \log_2\left(\frac{2\bar{M}}{M_0}\right)$.
Let $N(k)$ denote the number of inner line-search iterations up to iteration $k$ of $\SAHBA$. Then, 
\begin{align*}
N(k)&=\sum_{j=0}^{k}(i_{j}+1)\leq i_0+1 + \sum_{j=1}^{k}\left(\log_{2}\left(\frac{L_{j}}{L_{j-1}}\right)+2\right) 
\leq 2(k+1) + 2 \log_2\left(\frac{2\bar{M}}{M_0}\right),
\end{align*}
since $L_{k} \leq 2\bar{M}= 2\max\{M,M_0\}$ in the last step. Thus, on average, the inner loop ends after two trials.


\subsubsection{Bound for the Number of Outer Iterations}
Our goal now is to establish per-iteration decrease of the potential $F_{\mu}$ if the stopping condition does not yet hold, and derive from that a global iteration complexity bound for $\SAHBA$.
Let us fix iteration counter $k$. As said, the main assumption of this subsection is that the stopping criterion is not satisfied, i.e., either
$\|v^{k}\|_{x^{k}} \geq \Delta_{k}$ or $\|v^{k-1}\|_{x^{k-1}} \geq \Delta_{k-1}$. 
Without loss of generality, we assume that the first inequality holds, i.e., $\|v^{k}\|_{x^{k}} \geq \Delta_{k}$, and consider iteration $k$. Otherwise, if the second inequality holds, the same derivations can be made considering the iteration $k-1$ and using the second inequality $\|v^{k-1}\|_{x^{k-1}} \geq \Delta_{k-1}$. Thus, at the end of the $k$-th iteration 
\begin{equation}
\label{eq:SO_per_iter_proof_1}
\|v^{k}\|_{x^{k}} \geq \Delta_{k}=\sqrt{\frac{\eps}{12L_{k}\nu}}.
\end{equation}
Since the step-size $\alpha_k= \min\{1,\frac{1}{ 2\norm{v^{k}}_{x^{k}}}\} = \ct_{\mu,L_k}(x^{k})$ in \eqref{eq:SO_alpha_k} satisfies $\alpha_k \leq 1$ and $\alpha_k \norm{v^{k}}_{x^{k}} \leq 1/2$ (see \eqref{eq:SO_t_opt_def} and a remark after it), we can repeat the derivations of \eqref{eq:progressSOM}, 
changing  \eqref{eq:cubicestimate} to  \eqref{eq:SO_LS_1}.
In this way we obtain the following counterpart of \eqref{eq:progressSOM} with $t=\alpha_k$, $L=L_k$, $x=x^k$, $v_{\mu,L_k}(x^{k}) \eqdef v^k$: 
\begin{align}
F_{\mu}(x^{k+1})-F_{\mu}(x^{k})&\leq - \norm{v^{k}}^{3}_{x^{k}} \frac{L_k\alpha_k^2}{12}\left(3 - 2\alpha_k - \frac{12 \mu}{L_k\norm{v^{k}}_{x^{k}}} \right) \leq - \norm{v^{k}}^{3}_{x^{k}} \frac{L_k\alpha_k^2}{12}\left(1 - \frac{12 \mu}{L_k\norm{v^{k}}_{x^{k}}} \right)\label{eq:SO_per_iter_proof_2},
\end{align}
where in the last inequality we used that $\alpha_k \leq 1$ by construction. 
Substituting $\mu = \frac{\eps}{4\nu}$, and using \eqref{eq:SO_per_iter_proof_1}, we obtain
\begin{align*}
1 - \frac{12 \mu}{L_k\norm{v^{k}}_{x^{k}}} &= 1 - \frac{12 \eps}{4\nu L_k\norm{v^{k}}_{x^{k}}}  
\stackrel{\eqref{eq:SO_per_iter_proof_1}}{\geq} 1 - \frac{3 \eps}{\nu L_k\sqrt{\frac{\eps}{12L_{k}\nu}}} = 1 - \frac{6 \sqrt{\eps}}{ \sqrt{3\nu L_k}} \geq 1 - \frac{6 \sqrt{\eps}}{ \sqrt{ 3 \cdot 144 \nu \eps	}}  \geq \frac{1}{2},
\end{align*}
using that, by construction, $L_k =2^{i_k}M_k \geq \underline{L} = 144 \eps$ and that $\nu \geq 1$. 
Hence, from \eqref{eq:SO_per_iter_proof_2}, 
\begin{equation}
\label{eq:SO_per_iter_proof_3}
F_{\mu}(x^{k+1})-F_{\mu}(x^{k})\leq  - \norm{v^{k}}^{3}_{x^{k}} \frac{L_k\alpha_k^2}{24}.
\end{equation}
Substituting into \eqref{eq:SO_per_iter_proof_3} the two possible values of the step-size $\alpha_k$ in \eqref{eq:SO_alpha_k} gives
\begin{equation}
\label{eq:SO_per_iter_proof_8}
F_{\mu}(x^{k+1})-F_{\mu}(x^{k})\leq 
\left\{
\begin{array}{ll}
- \norm{v^{k}}^{3}_{x^{k}} \frac{L_k}{24}, & \text{if }  \alpha_k=1,\\
&\\
-  \norm{v^{k}}_{x^{k}} \frac{L_k}{96}, & \text{if } \alpha_k=\frac{1}{2\norm{v^{k}}_{x^{k}}}.
\end{array}\right.
\end{equation}
This implies
\begin{equation}
\label{eq:SO_per_iter_proof_7}
F_{\mu}(x^{k+1})-F_{\mu}(x^{k}) \leq - \frac{L_k\norm{v^{k}}_{x^{k}}}{96} \min\left\{1, 4\norm{v^{k}}_{x^{k}}^2 \right\} \eqdef -\delta_{k}.
\end{equation}
Rearranging and summing these inequalities for $k$ from $0$ to $K-1$, and using that $L_k \geq \underline{L}$, we obtain
\begin{align}
K\min_{k=0,...,K-1} &\frac{\underline{L}\norm{v^{k}}_{x^{k}}}{96} \min\left\{1, 4\norm{v^{k}}_{x^{k}}^2\right\} \leq  \sum_{k=0}^{K-1} \delta_k 
\leq F_{\mu}(x^{0})-F_{\mu}(x^{K}) \notag \\
&\stackrel{\eqref{eq:potential}}{=} f(x^{0}) - f(x^{K}) + \mu (h(x^{0}) - h(x^{K})) \leq f(x^{0}) - f_{\min}(\setX) + \eps, \label{eq:SO_per_iter_proof_9}
\end{align}
where we used that, by the assumptions of Theorem \ref{Th:SOAHBA_conv}, $x^{0}$ is a $4\nu$-analytic center defined in \eqref{eq:analytic_center} and $\mu = \frac{\eps}{4\nu}$, implying that $h(x^{0}) - h(x^{K}) \leq 4\nu = \eps/\mu$.
Thus, up to passing to a subsequence, we have $\norm{v^{k}}_{x^{k}} \to 0$ as $k \to \infty$, which makes the stopping criterion in Algorithm \ref{alg:SOAHBA} achievable.

Assume now that the stopping criterion does not hold for $K$ iterations of $\SAHBA$. 
Then, for all $k=0,\ldots,K-1,$ it holds that 
\begin{align}
\delta_k &= \frac{L_k}{96} \min\left\{\norm{v^{k}}_{x^{k}}, 4\norm{v^{k}}_{x^{k}}^3 \right\} 
\stackrel{\eqref{eq:SO_per_iter_proof_1}}{\geq} \frac{L_k}{96} \min\left\{  \sqrt{\frac{\eps}{12L_{k}\nu}} ,  \frac{4 \eps^{3/2}}{12^{3/2}L_{k}^{3/2}\nu^{3/2}}  \right\} \notag \\
&\stackrel{L_k \leq 2\bar{M}, \nu \geq 1}{\geq}  \frac{1}{96} \min\left\{  \frac{L_k \sqrt{\eps}}{\sqrt{24\bar{M}} \nu^{3/2}} ,  \frac{ \eps^{3/2}}{2 \cdot 3^{3/2} L_{k}^{1/2}\nu^{3/2}}  \right\}  \notag \\
&\stackrel{L_k \leq 2\bar{M},L_k \geq 144  \eps}{\geq} \frac{1}{96} \min\left\{ \frac{(144 \eps) \cdot \sqrt{\eps}}{\sqrt{24\bar{M}} \nu^{3/2}} ,  \frac{ \eps^{3/2}}{6\sqrt{6\bar{M}}\nu^{3/2}}  \right\}  = \frac{\eps^{3/2}}{576 \nu^{3/2} \sqrt{6\bar{M}} }.
\end{align}
Thus, from \eqref{eq:SO_per_iter_proof_9}
\begin{align*}
K \frac{\eps^{3/2}}{576 \nu^{3/2} \sqrt{6\bar{M}} } \leq f(x^{0}) - f_{\min}(\setX) + \eps. 
\end{align*}
Hence, recalling that $\bar{M}=\max\{M_0,M\}$, we obtain $K \leq \frac{576 \nu^{3/2} \sqrt{6\max\{M_0,M\}}(f(x^{0}) - f_{\min}(\setX)+ \eps)}{\eps^{3/2} }$.
Thus, we obtain the bound on the number of iterations on which the stopping criterion is not satisfied.
This, combined with the bound for the number of inner steps in Section \ref{sec:backtrack2}, proves the complexity bound statement of Theorem \ref{Th:SOAHBA_conv}.

\subsubsection{Generating $(\eps_{1},\eps_{2})$-2KKT Point}
To finish the proof of Theorem \ref{Th:SOAHBA_conv}, we show that if the stopping criterion in Algorithm \ref{alg:SOAHBA} holds, i.e., $\norm{v^{k-1}}_{x^{k-1}} < \Delta_{k-1}$ and $\norm{v^{k}}_{x^{k}} < \Delta_{k}$, then the algorithm has generated an $(\eps_{1},\eps_{2})$-2KKT point of \eqref{eq:Opt} according to Definition \ref{def:eps_SOKKT}, with $\eps_{1}=\eps$ and $\eps_{2}=\frac{\max\{M_{0},M\}\eps}{24\nu}$.

Clearly, \eqref{eq:eps_SO_optim_equality_set} in Definition \ref{def:eps_SOKKT} holds by the construction of the algorithm. Thus, we focus on showing \eqref{eq:eps_SO_optim_FO} and \eqref{eq:eps_SO_optim_SO}.

Let the stopping criterion hold at iteration $k$. First, we focus on showing \eqref{eq:eps_SO_optim_FO}. Using the first-order optimality condition \eqref{eq:opt1} for the subproblem \eqref{eq:SO_finder} solved at iteration $k-1$, there exists a Lagrange multiplier $y^{k-1}\in\R^{m}$ such that \eqref{eq:opt1} holds. Now, expanding the definition of the potential \eqref{eq:potential} and adding $\nabla f(x^{k})$ to both sides, we obtain from \eqref{eq:opt1} 
\begin{align*}
 \nabla f(x^{k}) - & \bA^{\ast}y^{k-1} + \mu \nabla h(x^{k-1}) \\
&=\nabla f(x^{k}) - \nabla f(x^{k-1}) - \nabla^2 f(x^{k-1})v^{k-1} - \frac{L_{k-1}}{2}\norm{v^{k-1}}_{x^{k-1}}H(x^{k-1})v^{k-1}.
\end{align*}
Setting $s^{k}\eqdef \nabla f(x^{k})-\bA^{\ast}y^{k-1} \in \setE^{\ast}$ and $g^{k-1}\eqdef-\mu\nabla h(x^{k-1})$, after multiplication by $[H(x^{k-1})]^{-1}$, this is equivalent to 
\begin{align*}
[H(x^{k-1})]^{-1}\left(s^{k}-g^{k-1}\right)&=[H(x^{k-1})]^{-1}\left(\nabla f(x^{k}) - \nabla f(x^{k-1}) - \nabla^2 f(x^{k-1})v^{k-1}- \frac{L_{k-1}}{2}\norm{v^{k-1}}_{x^{k-1}}H(x^{k-1})v^{k-1}\right).
\end{align*}
Multiplying both of the above equalities, we arrive at 
\begin{align*}
\left(\norm{s^{k}-g^{k-1}}^{\ast}_{x^{k-1}}\right)^{2}
&=\left( \left\|\nabla f(x^{k}) - \nabla f(x^{k-1}) - \nabla^2 f(x^{k-1})v^{k-1}  - \frac{L_{k-1}}{2}\norm{v^{k-1}}_{x^{k-1}}H(x^{k-1})v^{k-1} \right\|_{x^{k-1}}^*\right)^{2}.
\end{align*}
Taking the square root and applying the triangle inequality, we obtain 
\begin{align}
\norm{s^{k}-g^{k-1}}^{\ast}_{x^{k-1}}&\leq\norm{\nabla f(x^{k}) - \nabla f(x^{k-1}) - \nabla^2 f(x^{k-1})v^{k-1}}^{\ast}_{x^{k-1}}+\frac{L_{k-1}}{2}\norm{v^{k-1}}^{2}_{x^{k-1}} \notag \\
&\stackrel{\eqref{eq:SO_LS_2}}{\leq} \frac{L_{k-1}}{2}\norm{\alpha_{k-1}v^{k-1}}^{2}_{x_{k-1}}+\frac{L_{k-1}}{2}\norm{v^{k-1}}^{2}_{x^{k-1}}. \label{eq:SO_eps_KKT_proof_-1}
\end{align}
Since the stopping criterion holds, at iteration $k-1$ we have
\begin{align}
 \|v^{k-1}\|_{x^{k-1}} < \Delta_{k-1} = \sqrt{\frac{\eps}{12L_{k-1}\nu}} \leq \sqrt{\frac{\eps}{12 \cdot 144 \eps \nu}} < \frac{1}{2}\label{eq:SO_eps_KKT_proof_0},
\end{align}
where we used that, by construction, $L_{k-1} \geq \underline{L} = 144 \eps$ and that $\nu \geq 1$. Hence, by \eqref{eq:SO_alpha_k}, we have that $\alpha_{k-1}=1$ and $x^{k}=x^{k-1}+v^{k-1}$. 
This, in turn, implies that 
\begin{equation}
\label{eq:SO_eps_KKT_proof_1}
\norm{s^{k}-g^{k-1}}^{\ast}_{x^{k-1}} \stackrel{\eqref{eq:SO_eps_KKT_proof_-1}}{\leq}  
L_{k-1}\norm{v^{k-1}}^{2}_{x^{k-1}} \leq L_{k-1} \Delta_{k-1}^2 = \frac{\eps}{12\nu},
\end{equation}
where we used the stopping criterion 
\[
\|v^{k-1}\|_{x^{k-1}} < \Delta_{k-1}=\sqrt{\frac{\eps}{12L_{k-1}\nu}}
\]
Recalling that $s^{k}\eqdef \nabla f(x^{k})-\bA^{\ast}y^{k-1} \in \setE^{\ast}$ and $g^{k-1}\eqdef-\mu\nabla h(x^{k-1})$, we obtain from \eqref{eq:SO_eps_KKT_proof_1}
\begin{equation}
\label{eq:SO_eps_KKT_proof_11}
\norm{-\frac{1}{\mu}(\nabla f(x^{k})-\bA^{\ast}y^{k-1})-\nabla h(x^{k-1})}^{\ast}_{x^{k-1}} \leq \frac{\eps}{12\mu\nu}=\frac{1}{3},
\end{equation}
where the last equality uses that $\mu=\frac{\eps}{4\nu}$.
Applying \eqref{eq:approx_normal_cone} with $s=-\frac{1}{\mu}(\nabla f(x^{k})-\bA^{\ast}y^{k-1})$ and  $\xi=\frac{\eps}{12\mu\nu}=\frac{1}{3}$, we obtain
\begin{equation}
\label{eq:SO_eps_KKT_proof_12}
\inner{-\frac{1}{\mu}(\nabla f(x^{k})-\bA^{\ast}y^{k-1}),x-x_k}<\nu+\frac{\sqrt{\nu}+\xi}{1-\xi} \xi = \nu + \frac{\sqrt{\nu}+1/3}{2}\qquad \forall x \in \bar{\setK},
\end{equation}
whence,  
\begin{equation}
\label{eq:SO_eps_KKT_proof_13}
\inner{\nabla f(x^{k})-\bA^{\ast}y^{k-1},x-x_k}>-\mu\nu-\mu\frac{\sqrt{\nu}+1/3}{2} > -\eps \qquad \forall x \in \bar{\setK},
\end{equation}
where we used that $\mu=\frac{\eps}{4\nu}$ and that $\nu\geq 1$. Thus, we obtain that \eqref{eq:eps_SO_optim_FO} holds.\\
Finally, we show the second-order condition \eqref{eq:eps_SO_optim_SO}.
By inequality \eqref{eq:PD} for subproblem \eqref{eq:SO_finder} solved at iteration $k$, we obtain on $\setL_0$
\begin{align}
\nabla^2 f(x^{k})  &\succeq - \frac{L_{k}\norm{v^{k}}_{x^{k}}}{2} H(x^{k}) \succeq -\frac{L_{k} \Delta_k}{2} H(x^{k}) \notag \\
& =  - \frac{L_{k}}{2}\sqrt{\frac{\eps}{12L_{k}\nu}} H(x^{k}) = - \frac{\sqrt{L_k \eps}}{(48\nu)^{1/2}}H(x^{k})  
\succeq - \frac{ \sqrt{2\bar{M}\eps}}{(48\nu)^{1/2}}H(x^{k})=- \frac{ \sqrt{\bar{M}\eps}}{(24\nu)^{1/2}}H(x^{k}), \label{eq:SO_eps_KKT_proof_4}
\end{align}
where we used the second part of the stopping criterion, i.e., $\norm{v^{k}}_{x^{k}}< \Delta_k$ and that $L_k \leq 2\bar{M}=2\max\{M,M_0\}$ (see Section \ref{sec:backtrack2}). Thus, \eqref{eq:eps_SO_optim_SO} holds with $\eps_2=\frac{\max\{M,M_0\}\eps}{24\nu}$, which finishes the proof of Theorem \ref{Th:SOAHBA_conv}.

\end{document}